\documentclass[english,10pt]{article}
\usepackage{amssymb}
\usepackage{tikz}
\usepackage{cancel}

\usepackage{comment}
\usepackage{enumitem}

\usepackage{lineno}  

\usepackage{ulem}

\def\0{\emptyset}

\def\p{\noindent{\bf Proof.~~}}
\def\q{\hfill\rule{1ex}{1ex}}

\newtheorem{Theorem}{Theorem}[section]

\newtheorem{Definition}[Theorem]{Definition}
\newtheorem{Example}[Theorem]{Example}
\newtheorem{Lemma}[Theorem]{Lemma}

\newtheorem{Observation}[Theorem]{Observation}
\newtheorem{Proposition}[Theorem]{Proposition}

\newcommand{\G}{\mbox{$\cal G$}}

\usepackage{algorithm,algcompatible,amsmath}
\usepackage{algpseudocode}
\algnewcommand\INPUT{\item[\textbf{Input:}]}%
\algnewcommand\OUTPUT{\item[\textbf{Output:}]}%

\usepackage{color}
\usepackage{CJK}

\begin{document}

\title{Optimal connectivity of second order iterated line graphs}

\author{Run Zou\thanks{College of Mathematics and System Sciences,
Xinjiang University, Urumqi, Xinjiang 830046, PRC. Email: zourun163163@163.com},  Wei Xiong\thanks{College of Mathematics and System Sciences,
Xinjiang University, Urumqi, Xinjiang 830046, PRC. Email: xingheng-1985@163.com},
Mingquan Zhan\thanks{Department of Mathematics,
Millersville University, Millersville, PA 17551, USA.
Email: Mingquan.Zhan@millersville.edu},
Hong-Jian Lai\thanks{School of Mathematics and Systems Science, Guangdong Polytechnic Normal University, Guangzhou 510665,
China, and Department of Mathematics, West Virginia University,
Morgantown, WV 26506, USA. Email: hjlai2015@hotmail.com}
}
\date{}
\maketitle

\begin{abstract}
The line graph $L(G)$ of a graph $G$ is defined to be the simple graph whose vertices
are the edges of $G$, where two vertices in $L(G)$ are adjacent if and only if the corresponding edges in $G$
are incident with a common vertex, and define $L^2(G)=L(L(G))$.
For positive integers $d$ and $k$, the function
 $\kappa_{L^2}(d,k) = \inf\{\kappa(L^2(G)): \kappa'(G) \ge k \mbox{ and } \delta(G) \ge d\}$
has been investigated.
Niepel and Knor proved that $\kappa_{L^2}(d,1)\geq d-1$, for any integer $d \ge 3$.
In this research, it is proved that if $d\geq 3$ and $k\geq 1$,
then $\kappa_{L^2}(d,k)= \min\{f(d,k), 4d-6\}$,
where
\begin{equation}
f(d,k) = \left\{
\begin{array}{ll}
k(d-k), & \mbox{ if $1\leq k\leq \lfloor\frac{d}{2}\rfloor$, }
\\
kd-k^2+2k(\lceil \frac{d}{2}\rceil)-(\lceil\frac{d}{2}\rceil)d,
& \mbox{ if $\lfloor\frac{d}{2}\rfloor< k <\frac{3d-1}{4}$, }
\\
kd-k^2+2k\lfloor \frac{d}{2}\rfloor-2(\lfloor \frac{d}{2}\rfloor)^2,
& \mbox{ if $\frac{3d-1}{4}\leq k<d$, }\\
d(\lceil \frac{d}{2}\rceil),
& \mbox{ if } d=k.
\end{array} \right.\nonumber
\end{equation}
\end{abstract}

{\small
}

\vskip 0.35cm

\section{Introduction}
Networks of various kinds are often modeled as graphs. Line graphs and iterated line graphs play useful roles in network modeling and analysis, as demonstrated in
 \cite{CLWJ22, DixD25+}, among others. Graph connectivity is one of the most fundamental and extensively studied invariants in graph theory, with numerous applications in network reliability and design. The study of connectivity in iterated line graphs has attracted considerable attention.
Chartrand and Stewart conducted pioneering work on the relationships between the edge and vertex connectivities of a graph
$G$ and its line graph $L(G)$ in \cite{ChSt69}. Knor and Niepel \cite{KnNi03} examined  the the connection  between  $L^2(G)$, the second order iterated line graph of
$G$, and the minimum degree of $G$. Several studies by Shao \cite{Shao05, Shao10, Shao18, Shao22+}
further investigated the connectivity properties of iterated line graphs and essential edge cuts in the original graph.
Additional research on the connectivity and related properties of iterated line graphs
can be found in \cite{BoKN06, FeSu++, Zamf70}.
The purpose of this work is to provide a more precise characterization of the connectivity of the second iterated line graph, thereby extending the results of
 Knor and Niepel \cite{KnNi03}.

In this research, we consider simple and finite graphs and follow the terminology and notation of  \cite{BoMu08}
for undefined terms.  In particular, we use $\kappa(G)$, $\kappa'(G)$, $\delta(G)$,  and $\Delta(G)$ to denote the connectivity, the edge-connectivity,
the minimum degree,  and the maximum degree of a graph $G$, respectively.
The {\bf line graph} $L(G)$ of a graph $G$ is defined to be the simple graph whose vertices
are the edges of $G$, where two vertices in $L(G)$ are adjacent if and only if the corresponding edges in $G$
are incident with a common vertex. A survey by Hemminger and Beineke   \cite{HeBe78} provides an excellent overview of general properties of line graphs.
 We define $L^0(G) = G$. For an integer $k \ge 1$, the $k$-th order of {\bf iterated line graph} of a graph $G$,
denoted $L^k(G)$, is defined recursively by $L^{k}(G)=L(L^{k-1}(G))$.

Let $G, H$ be  graphs. Define $G \cup H$ to be the graph with vertex set $V(G \cup H) = V(G) \cup V(H)$
and edge set $E(G) \cup E(H)$. If $E' \subseteq E(G)$, then we define $V(E')$ to be the set of vertices of $G$ that
are incident with at least one edge in $E'$. When $E' = \{e\}$, we often use $V(e)$ for $V(\{e\})$.
The subgraph of $G$ induced by $E'$, denoted $G[E']$, has vertex set
$V(E')$ and edge set $E'$.

Let $G$ be a connected graph. For vertex subsets $U,V\subseteq V(G)$,   define $E_G[U,V] = \{uv \in E(G): u \in U, v \in V\}$,
and let $\partial_G(U) = E_G[U, V(G)-U]$. If $U = \{v\}$, then we often write  $\partial_G(v)$ and $E_G[v,V]$ for $\partial_G(\{v\})$ and $E_G[\{v\},V]$, respectively.
A graph is {\bf trivial} if it contains no edges. An edge cut $X$ of $G$ is a subset of edges of the form $X = E_G[U, V(G) - U]$ for some vertex subset $U \subseteq V(G)$.
By definition, if $G$ is a connected graph, then an edge cut
\begin{equation} \label{min-edge-cut}
\mbox{$X = E_G[U, V(G) - U]$ is minimal if and only if both $G[U]$ and $G[V(G) - U]$ are connected.}
\end{equation}
Suppose $G$ is a connected graph with $|E(G)|\geq 2$. An edge cut $X$ is {\bf essential}  if $G - X$ has at least two nontrivial components.
For an integer $k>0$, a graph $G$ is {\bf essentially $k$-edge-connected} if $G$ does not have an essential edge cut of size less than $k$.
The {\bf essential edge connectivity} of a graph $G$, denoted by $ess'(G)$, is the maximum value of $k$ for
which $G$ is essentially $k$-edge connected. 
Observation \ref{obs-1} follows directly from these  definitions.

\begin{Observation} \label{obs-1}
(Proposition 1.1.3 of \cite{Shao05})
Let $G$ be a connected graph. Then $\kappa(L(G))\geq k$ if and only if $ess'(G)\geq k$.
\end{Observation}

By the definition of a line graph, if $G$ is a path, a cycle, or a $K_{1,3}$, then the line graph is also a path or a cycle.
Observation \ref{obs-1a} below follows immediately from definitions.

\begin{Observation} \label{obs-1a}
Let $G$ be a connected graph which is not isomorphic to a path, or a cycle or a $K_{1,3}$. Then $L(G)$ not isomorphic to a path, or a cycle or a $K_{1,3}$.
\end{Observation}

Therefore, throughout this paper, we always assume that $G$ is not isomorphic to a path, a cycle, or a $K_{1,3}$.

Let $G$ be a connected simple graph with minimum degree $\delta$. A number of studies investigated
how the minimum degree, the connectivity and the edge-connectivity of a graph $G$ can be used
to predict or describe the connectivity properties of $L^2(G)$.
In 1969, Chartrand and Stewart \cite{ChSt69} studied the edge and vertex connectivity relations between a graph G and its line graph $L(G)$.
Knor and Niepel \cite{KnNi03} showed that
if G is a connected graph with $\delta \ge 3$, then $L^2(G)$ is $(\delta-1)$-connected
and showed $\kappa(G)\geq 4$, then $\kappa(L^2(G))\geq 4d-6$.
More recent studies on the connectivity of iterated line graphs
were done by  Knor and Niepel \cite{KnNi06}. Shao in \cite{Shao10} used essential edge
connectivity of G to study the connectivity of a iterated line graph.

\begin{Theorem} (Knor and Niepel \cite{KnNi03})\label{Knor}
Let G be a connected graph with $\delta \ge 3$, then $L^2(G)$ is $(\delta-1)$-connected.
\end{Theorem}

Theorem \ref{Knor} motivates the construction of the following examples. These examples are related to the discussions
on the connectedness of second order iterated line graphs.

\begin{center}
  \begin{tikzpicture} [scale=0.5]

   \path    (3,0)  coordinate (v);  \draw [fill=black] (v) circle (0.08cm); \node [below] at (v) { $v_0$};

 \draw (0,0) ellipse (1.5cm and 3.5cm);      \draw (6,0) ellipse (1.5cm and 3.5cm);

        \foreach \i in {-4,...,4}
    {   \path    (0,0.68*\i)  coordinate (v\i);  \draw [fill=black] (v\i) circle (0.08cm);
     \path    (6,0.68*\i)  coordinate (u\i);  \draw [fill=black] (u\i) circle (0.08cm);
}
    \node [left] at (v0) { $v_1^1$};

        \foreach \i in {0}
    { \draw (v)--(v\i);
}

        \foreach \i in {-3,...,3}
    { \draw (v)--(u\i);
}

 \node [right] at (5.4,0) { $\left.\begin{array}{c}  \\
  \\
  \\ \\
   \end{array}\right\} d-1$};

\node [below] at (0,-3.5) {$J_1=K_{n_1}$}; \node [below] at (6,-3.5) {$J_2=K_{n_2}$};
\node [below] at (3,-5) {$J(d,1,n_1,n_2)$};

     \end{tikzpicture} \hskip1cm
       \begin{tikzpicture} [scale=0.5]

   \path    (3,0)  coordinate (v);  \draw [fill=black] (v) circle (0.08cm); \node [below] at (v) { $v_0$};

 \draw (0,0) ellipse (1.5cm and 3.5cm);       \draw (6,0) ellipse (1.5cm and 3.5cm);

         \foreach \i in {-4,...,4}
    {   \path    (0,0.68*\i)  coordinate (v\i);  \draw [fill=black] (v\i) circle (0.08cm);
     \path    (6,0.68*\i)  coordinate (u\i);  \draw [fill=black] (u\i) circle (0.08cm);
}

            \foreach \i in {-1,...,1}
    { \draw (v)--(v\i);
}

        \foreach \i in {-2,...,2}
    { \draw (v)--(u\i);
}

 \node [right] at (5.4,0) { $\left.\begin{array}{c}  \\
  \\
  \\
   \end{array}\right\} d-k$};

 \node [left] at (0.6,0) { $k \left\{\begin{array}{c}  \\
   \end{array}\right. $};

\node [below] at (0,-3.5) {$J_1=K_{n_1}$}; \node [below] at (6,-3.5) {$J_2=K_{n_2}$};
\node [below] at (3,-5) {$J(d,k,n_1,n_2)$};

     \end{tikzpicture}

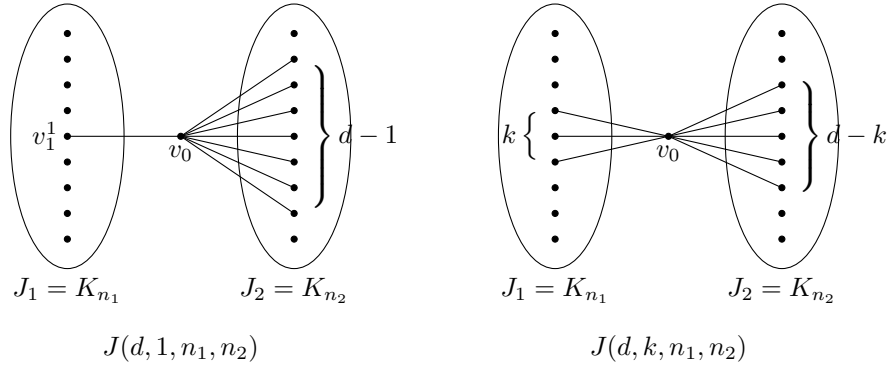
     Figure 1: Graphs in Example 1.4.
\end{center}

\begin{Example} \label{ex-1}
Let $d, k, n_1, n_2$ be integers with $\min\{n_1, n_2\} \ge d+1 \ge 4$ and $d > k >0$, and for $i \in \{1,2\}$, $J_i \cong K_{n_i}$ be a complete graph with $V(J_i) = \{v_1^i, v_2^i, ..., v_{n_i}^i\}$. Define $J(d,k,n_1, n_2)$ to be a graph obtained from the disjoint union $J_1 \cup J_2$ by adding a new vertex $v_0$ and new edges $\{v_0v_j^1: 1 \le j \le k\} \cup \{v_0v_j^2: 1 \le j \le d-k\}$. (See Figure 1 for illustrations).
 
Let $G = J(d,k, n_1, n_2)$. Denote $E_1 = \{v_0v_j^1: 1 \le j \le k\}$ and $E_2 = \{v_0v_j^2: 1 \le j \le d-k\}$. Then $E_1 \cup E_2 \subseteq V(L(G))$. As $ J(d,k,n_1, n_2) - \{v_0\}$ is a disjoint union of complete graphs $K_{n_1}$ and $K_{n_2}$,the edge set $W = \{e_1e_2: e_1 \in E_1, e_2 \in E_2\}$ in $L (G)$ is an essential edge cut of $L(G)$ with $|W| = k(d-k)$. It follows that  $\kappa(L^2(G)) \leq k(d-k)$. In particular, when $k=1$, we have $\kappa(L^2(G)) \leq d-1$.
\end{Example}

For each ordered pair of integers $(d,k)$ with $d\ge k\ge 1$ and  $d \ge 3$, define
\begin{equation} \label{def-g-k}
\begin{array}{ll}
\G(d,k)  & = \{G: \text{$G$ is a simple graph with $\kappa'(G)=k\geq 1$ and $\delta(G)= d\geq 3$}\}
\\
\kappa_{L^2}(d,k)  & = \inf\{ \kappa(L^2(G)): G\in \G(d,k)\}.
\end{array}
\end{equation}
Thus Theorem \ref{Knor} suggests that, if $d \ge 3$, then $\kappa_{L^2}(d,1) \ge d-1$.
This, together with Example \ref{ex-1}, implies that $\kappa_{L^2}(d,1)= d-1$.
This motivates the current research. The purpose of this research aims to extend Theorem \ref{Knor} by Knor and Niepel
to graphs with all integers $d \ge 3$ and $k\geq 1$. The following theorem is the main result of this research.
\begin{Theorem} \label{main-1}
Let $d$ and $k$ be integers with $d\ge k\ge 1$ and $d \ge 3$. Then
\begin{equation}\label{function}
\kappa_{L^2}(d,k)=\min\{f(d,k),4d-6\},\nonumber
\end{equation}
where
\begin{equation}
f(d,k) = \left\{
\begin{array}{ll}
k(d-k),  & \mbox{ if $1\leq k\leq \lfloor\frac{d}{2}\rfloor$, }
\\
kd-k^2+2k(\lceil \frac{d}{2}\rceil)-(\lceil\frac{d}{2}\rceil)d,
& \mbox{ if $\lfloor\frac{d}{2}\rfloor< k <\frac{3d-1}{4}$, }
\\
kd-k^2+2k\lfloor \frac{d}{2}\rfloor-2(\lfloor \frac{d}{2}\rfloor)^2,
& \mbox{ if $\frac{3d-1}{4}\leq k<d$, }\\
d(\lceil \frac{d}{2}\rceil),
& \mbox{ if d=k. }
\end{array} \right.\nonumber
\end{equation}
\end{Theorem}
In the next section, we present the preliminaries related to this problem.  The proof of the main result will be given in the last section.


\section{Essential edge cuts in $L(G)$}
In \cite{Shao10}, for a given essential edge cut $X$ of a line graph $L(G)$, Shao introduced the vertex subsets $V_1(X), V_2(X)$ and $V_{12}(X)$ of $G$
as in Definition \ref{v-x} below.

\begin{Definition} \label{v-x}
Let $G$ be a connected graph that is not isomorphic to a cycle, a path or a $K_{1,3}$.
\\
(i) For each edge $e\in E(G)$, let $v_e\in V(L(G))$ denote the vertex corresponding to $e$ in the process of forming $L(G)$ from $G$.
\\
(ii) Let $X$ be a minimal
essential  edge cut of $L(G)$, and let $L_1,L_2$ be the two components of $L(G)-X$.
Define $f_X :E(G) \to \{1,2\}$ be the ({\bf $X$-induced}) 2-edge-coloring of $G$ such that, for each $i \in \{1,2\}$,  $f_X(e)=i$ if and only if $v_e\in V(L_i)$.
For each vertex subset $W\subseteq V(G)$, define
\[
E_{i, X}(W)=\{e\in \partial_G(W): f_X(e)=i\} \mbox{ and } N_{i,X}(W) = \{v \in V(G) - W: E_G[v, W] \subseteq E_{i, X}(W)\}.
\]
When $W=\{v\}$, we often use $E_{i,X}(v)$ and $N_{i,X}(v)$ for $E_{i,X}(\{v\})$ and $N_{i,X}(\{v\})$, respectively.
\\
(iii) A vertex $v\in V(G)$ is \textbf{$X$-mono-colored} if for some $i \in \{1,2\}$, $\partial_G(v) = E_{i,X}(v)$. Define, for $i \in \{1,2\}$,
\begin{equation}\label{V12}
\begin{array}{ll}
V_{i}(X)  & = \{v\in V(G): \partial_G(v) = E_{i,X}(v) \},
\\
V_{12}(X) & = \{v\in V(G): \partial_G(v) \cap E_{1,X}(v) \neq \emptyset \mbox{ and } \partial_G(v) \cap E_{2,X}(v) \neq \emptyset \}.
\end{array}
\end{equation}
\end{Definition}
When the essential edge cut $X$ is understood from the context, we often omit the subscript $X$
in the notation above, using $f$ for $f_X$, $E_{i}$ for $E_{i,X}$ and saying that a vertex is mono-colored in stead of being $X$-mono-colored.

As an example, for the graphs $G$ and $L(G)$ in Figure 2, $X=\{e_0e_7, e_4e_7, e_5e_8, e_6e_9\}$ is a minimal essential edge cut of $L(G)$, and the components of $L(G)-X$ satisfy
$V(L_1)=\{e_0, e_1, e_2, \cdots, e_6\}$ and $V(L_2)=\{e_7,e_8,\cdots,e_{12}\}$.  Thus $b_0$ in $G$ is  $X$-mono-colored since $\partial_G(b_0)=\{e_1,e_3,e_4\}= E_{1,X}(b_0)$. Also
 $b_2$ in $G$ is  also  $X$-mono-colored since  $\partial_G(b_2)=\{e_7,e_{11},e_{12}\}= E_{1,X}(b_2)$. In addition,
$V_1(X)=\{a_0,b_0,c_0\}$, $V_2(X)=\{a_2,b_2,c_2\},$ and  $X_{12}(X)=\{v_1,v_2,v_3\}$.
Note that as $G$ is connected, the set $V_{12}(X)\neq \emptyset$.
Observation \ref{Obs(G)} below follows from the definitions.

    \begin{center}
\begin{tikzpicture} [scale=0.7]
\foreach \i in {0,1,2}
{  \path    (2*\i,0)  coordinate (v\i);  \draw [fill=black] (v\i) circle (0.05cm);
  \path    (2*\i,1.5)  coordinate (a\i);  \draw [fill=black] (a\i) circle (0.05cm);
    \path    (2*\i,3)  coordinate (b\i);  \draw [fill=black] (b\i) circle (0.05cm);
    }
    \draw (v0)--(v2) (a0)--(a2) (b0)--(b2);
    \draw (v0)--(b0) (v0).. controls (-1, 0.5) and (-1, 2.5) ..(b0);

 \node [below] at (v0) {$c_0$};  \node [below] at (v1) {$v_1$};     \node [below] at (v2) {$c_2$};
  \node [below] at (a0) {$a_0$};  \node [below] at (a1) {$v_2$};  \node [below] at (a2) {$a_2$};
   \node [above] at (b0) {$b_0$};   \node [above] at (b1) {$v_3$};  \node [above] at (b2) {$b_2$};

    \draw (v2)--(b2) (v2).. controls (5, 0.5) and (5, 2.5) ..(b2);

    \node [above] at (1,-0.2) { $e_6$};      \node [above] at (3,-0.2) { $e_9$};
    \node [above] at (1,1.3) {$e_5$};      \node [above] at (3,1.3) { $e_8$};
    \node [above] at (1,2.8) {$e_4$};      \node [above] at (3,2.8) { $e_7$};

    \node [right] at (-0.1, 0.75) { $e_2$};      \node [right] at (-0.1, 2.25) { $e_3$};
 \node [left] at (-0.6, 1.5) {$e_1$};

        \node [left] at (4.1, 0.75) { $e_{10}$};      \node [left] at (4.1, 2.25) {$e_{11}$};
 \node [right] at (4.5, 1.5) {$e_{12}$};
 \draw (a0)--(b1);  \node [right] at (1, 2.5) {$e_{0}$};
    \node [below] at (2,-1) {$G$};

     \end{tikzpicture} \hskip2cm
     \begin{tikzpicture} [scale=0.7]
\foreach \i in {0,1,2,3}
{  \path    (2*\i,0)  coordinate (v\i);  \draw [fill=black] (v\i) circle (0.05cm);
  \path    (2*\i,1.5)  coordinate (a\i);  \draw [fill=black] (a\i) circle (0.05cm);
    \path    (2*\i,3)  coordinate (b\i);  \draw [fill=black] (b\i) circle (0.05cm);
    }

\draw (v1)--(v2) (a1)--(a2) (b1)--(b2);
\draw (v0)--(b0)  (v0).. controls (-1, 0.5) and (-1, 2.5) ..(b0);

\draw (v3)--(b3)  (v3).. controls (7, 0.5) and (7, 2.5) ..(b3);

      \node [below] at (2, 0) { $e_6$};       \node [below] at (4, 0) { $e_9$};
       \node [below] at (2, 1.5) {$e_5$};       \node [below] at (4, 1.5) { $e_8$};
            \node [above] at (2, 3) {$e_4$};       \node [above] at (4, 3) { $e_7$};

          \node [below] at (0, 0) {$e_2$};       \node [left] at (0, 1.5) {$e_1$};
            \node [above] at (0, 3) { $e_3$};
            \draw (v1)--(v0) (v1)--(a0);           \draw (a1)--(v0) (a1)--(b0);
               \draw (b1)--(a0) (b1)--(b0);

           \node [below] at (6, 0) { $e_{10}$};       \node [right] at (6, 1.5) {$e_{12}$};
            \node [above] at (6, 3) { $e_{11}$};

                \draw (v2)--(v3) (v2)--(a3);                    \draw (a2)--(b3) (a2)--(v3);
                            \draw (b2)--(b3) (b2)--(a3);
    \node [below] at (3.5,-1) {$L(G)$ and $X=\{e_0e_7, e_4e_7, e_5e_8, e_6e_9\}$};

    \path    (2, 2.25)  coordinate (b);  \draw [fill=black] (b) circle (0.05cm);    \node [right] at (b) { $e_0$};
    \draw (b)--(b1) (b)--(b2) (b)--(v0) (b)--(b0) (b)--(a1);

     \end{tikzpicture}

     Figure 2: Illustration of Definition \ref{v-x}
\end{center}

\begin{Observation}\label{Obs(G)}
Let $G$ be a connected graph,
$X$ be a minimal essential  edge cut of $L(G)$. For each edge $e\in E(G)$, let $v_e$ denote the
vertex in $L(G)$ corresponding to $e$.
Each of the following holds:\\
(i) (Proposition 2.2(ii) of \cite{Shao10}) $v_{e_1}v_{e_2}\in X$ if and only if $e_1,e_2$
incident with common vertex in $V_{12}(X)$ and $f(e_1)\neq f(e_2)$.
\\
(ii) (Proposition 2.3(ii) of \cite{Shao10}) Each vertex of $G-V_{12}(X)$ is mono-colored in $G$, and moreover,
for each component $H$ of $G-V_{12}(X)$, all edges with at least one
end in $H$ have the same color as the edges of $H$.
\\
(iii) For each $v\in V_{12}(X)$, $d_G(v)=|E_{1,X}(v)|+|E_{2,X}(v)|$.
\\
(iv)  If for $i=1,2$, $V_i(X)\neq \emptyset$, then $E_i(V_{12}(X))$ 
contains an edge cut $E_G[V_i(X),V(G)-V_i(X)]$ of $G$.
Moreover, if $X$ is a minimal essential edge cut of $L(G)$, then $V_1(X) \cup V_2(X) \neq \emptyset$.
\\
(v) If $V_i(X) = \emptyset$ for some $i \in \{1,2\}$, then $E(G[V_{12}(X)]) \neq \emptyset$.
\end{Observation}

\p As $X$ is a minimal essential  edge cut of $L(G)$, in Definition \ref{v-x}(ii), $f_X$ is a 2-coloring of $E(G)$ with
$f_X^{-1}(1) \neq \emptyset$ and  $f_X^{-1}(2) \neq \emptyset$. Thus (iii) follows.
As for each $i \in \{1,2\}$, every path in $G$ linking a vertex in $V_1(X)$ to a vertex in $V_2(X)$ must us
an edge in $E_i(V_{12}(X))$, $E_i(V_{12}(X))$ contains an edge cut $E_G[V_i(X),V(G)-V_i(X)]$ of $G$.
If $V_i(X) = \emptyset$, then $E(G[V_{12}(X)] \subseteq f^{-1}(i) \neq \emptyset$, and so (iv) holds.
\q

By Observation \ref{Obs(G)}(i), if $X$ is an essential edge cut of $L(G)$, then
\begin{equation}\label{eq}
|X|=\sum \limits_{v\in V_{12}(X)}|E_{1,X}(v)| \cdot |E_{2,X}(v)| \ge 1.
\end{equation}
For example,  let $G = G(d,k,n_1, n_2)$ be a graph defined in Example \ref{ex-1}. Using the notation in Example \ref{ex-1}, we observe that, by (\ref{eq}),
$L(G)$ has a minimal essential edge cut $X = \{v_ev_{e'}: e \in E_G[v, V(J_1)], e' \in E_G[v, V(J_2)]\}$
with $|X|=k(d-k)$.  Define
\begin{equation} \label{def-g0}
\mbox{ ${\cal G}_0(d,k)=\{G\in {\cal G}(d,k): \kappa(L^2(G))= \kappa_{L^2}(d,k)\}$.}
\end{equation}
Since for each $G\in {\cal G}(d,k)$, $\kappa(L^2(G))$ is an integer, we conclude that for each integral pair $(d,k)$,
the value $\inf\{ \kappa(L^2(G)): G\in \G(d,k)\}$ is reachable by some graphs in $\G(d,k)$, and so ${\cal G}_0(d,k) \neq \emptyset$. Example \ref{ex-1} leads to the following.
\begin{Proposition} \label{kd-k2}
If $G\in \G_0(d,k)$, then $\kappa (L^2(G))\leq k(d-k).$
\end{Proposition}

\begin{Proposition}\label{m}
Let $G\in \G_0(d,k)$ be a graph. 
If $X$ is a minimum essential edge cut of $L(G)$
with $|V_{12}(X)|=m$, then $1\leq m \leq k$.
\end{Proposition}

\p Denote $V_{12}(X) = \{v_1, v_2, \ldots, v_m\}$. By \eqref{eq}, as $X$ is an essential edge cut, $m \ge 1$. If $m> k$, then
\[
\begin{split}
|X|&=\sum \limits_{v\in V_{12}(X)}|E_{1,X}(v)||E_{2,X}(v)| =\sum \limits_{v\in V_{12}(X)}|E_{1,X}(v)|(d_{G}(v)-|E_{1,X}(v)|)\\
&\geq \sum \limits_{v\in V_{12}(X)}1\cdot(d_{G}(v)-1) \geq \sum \limits_{v\in V_{12}(X)}(d-1)=|V_{12}(X)|(d-1)\\
&=m(d-1)>k(d-1)\geq k(d-k),
\end{split}
\]
contrary to Proposition \ref{kd-k2}.
\q

As suggested by  Proposition \ref{m}, in the rest of the discussions, we have the following definition.
\begin{Definition} \label{BS}
Let $X=\{v_1,v_2,\cdots,v_m\}$ be  an  essential edge cut of $L(G)$ with $|V_{12}(X)|=m$ such that $1\leq m< k$.
For each $v_i \in V_{12}(X)$, define
\begin{eqnarray*}
A_{i,1,X} & = & \left\{
\begin{array}{ll}
E_{1,X}(v_i), & \mbox{ if } |E_{1,X}(v_i)|\leq |E_{2,X}(v_i)|
\\
E_{2,X}(v_i), & \mbox{ if } |E_{1,X}(v_i)|>|E_{2,X}(v_i)|,
\end{array} \right.
\\
A_{i,2,X} & = & \left\{
\begin{array}{ll}
E_{1,X}(v_i), & \mbox{ if } |E_{1,X}(v_i)|> |E_{2,X}(v_i)|
\\
E_{2,X}(v_i), & \mbox{ if  } |E_{1,X}(v_i)|\leq |E_{2,X}(v_i)|,
\end{array} \right.
\end{eqnarray*}
$ a_{i,1,X}=|A_{i,1,X}|,$ $a_{i,2,X}=|A_{i,2,X}|$, and
\begin{eqnarray*}
N_{i,1,X} & = & \left\{
\begin{array}{ll}
N_{1,X}(v_i), & \mbox{ if } |E_{1,X}(v_i)|\leq |E_{2,X}(v_i)|
\\
N_{2,X}(v_i), & \mbox{ if } |E_{1,X}(v_i)|>|E_{2,X}(v_i)|,
\end{array} \right.
\\
N_{i,2,X} & = & \left\{
\begin{array}{ll}
N_{1,X}(v_i), & \mbox{ if } |E_{1,X}(v_i)|> |E_{2,X}(v_i)|
\\
N_{2,X}(v_i), & \mbox{ if  } |E_{1,X}(v_i)|\leq |E_{2,X}(v_i)|.
\end{array} \right.
\end{eqnarray*}
\end{Definition}
For instance,  for the graph $G$ and $L(G)$ shown in Figure 3,
let $X=\{e_1e_2,e_3e_4,e_3e_5, e_3e_6,e_3e_{10}, e_8e_7, e_8e_9,$ $e_{10}e_7,e_{10}e_9\}$ be a
minimal essential edge cut in $L(G)$. Then $V_{12}(X)=\{v_1,v_2,v_3\}$,
$A_{1,1,X}=E_{1,X}(v_1)=\{e_1\},
 A_{2,1,X}=E_{1,X}(v_2)=\{e_3\}, A_{3,1,X}=E_{1,X}(v_3)=\{e_7,e_9\}$, and
$A_{1,2,X}=E_{2,X}(v_1)=\{e_2\}, A_{2,2,X}=E_{2,X}(v_2)=\{e_4, e_5,e_6, e_{10}\}, A_{3,2,X}=E_{2,X}(v_3)=\{e_8,e_{10}\}$.

\begin{center}
\begin{tikzpicture} [scale=0.5]

\draw (0,2.3) ellipse (0.5cm and 3cm); \draw (4,2.3) ellipse (0.5cm and 3cm);   \node [above left] at (2.2,1.8) {\tiny $v_2$};   \node [above right] at (1.7,3.8) {\tiny $v_3$};
\node [above] at (2,0) {\tiny $v_1$};

 \path    (2,0)  coordinate (v1);  \draw [fill=black] (v1) circle (0.07cm);  \path    (0,0)  coordinate (a1);  \draw [fill=black] (a1) circle (0.07cm);
  \path    (4,0)  coordinate (b1);  \draw [fill=black] (b1) circle (0.07cm); \draw (a1)--(b1);

  \path    (2,2)  coordinate (vi);  \draw [fill=black] (vi) circle (0.05cm);
    \path    (0,2)  coordinate (a2);  \draw [fill=black] (a2) circle (0.05cm); \draw (vi)--(a2);
      \path    (4,2)  coordinate (b2);  \draw [fill=black] (b2) circle (0.05cm); \draw (vi)--(b2);
        \path    (4,2.5)  coordinate (b21);  \draw [fill=black] (b21) circle (0.05cm); \draw (vi)--(b21);
        \path    (4,1.5)  coordinate (b22);  \draw [fill=black] (b22) circle (0.05cm); \draw (vi)--(b22);

   \path    (2,4)  coordinate (vj);  \draw [fill=black] (vj) circle (0.05cm);\draw (vi)--(vj);
      \path    (0,4)  coordinate (a3);  \draw [fill=black] (a3) circle (0.05cm);\draw (vj)--(a3);
      \path    (4,4)  coordinate (b3);  \draw [fill=black] (b3) circle (0.05cm);\draw (vj)--(b3);
      \path    (0,4.5)  coordinate (a31);  \draw [fill=black] (a31) circle (0.05cm);\draw (vj)--(a31);

    \node [above] at (1,-0.2) {\tiny $e_1$};      \node [above] at (3,-0.2) {\tiny $e_2$};

    \node [above] at (1,1.3) {\tiny $e_3$};          \node [above] at (3,1.3) {\tiny $e_4$};    \node [above] at (3.5,1.7) {\tiny $e_5$};    \node [above] at (3,2) {\tiny $e_6$};

    \node [above] at (1,3.5) {\tiny $e_9$};      \node [above] at (3,3.8) {\tiny $e_8$};     \node [above] at (1,4) {\tiny $e_7$};
      \node [right] at (1.8,3) {\tiny $e_{10}$};

                                                         \begin{scope}[shift={(a1)}]
                    \draw (a1)--(100:0.5cm);          \draw (a1)--(80:0.5cm);
                                          \end{scope}

                                                         \begin{scope}[shift={(b1)}]
                    \draw (b1)--(100:0.5cm);          \draw (b1)--(80:0.5cm);
                                          \end{scope}

                                                                                                   \begin{scope}[shift={(a2)}]
                    \draw (a2)--(170:0.5cm);          \draw (a2)--(190:0.5cm);
                                          \end{scope}

                                                                                                             \begin{scope}[shift={(b2)}]
                    \draw (b2)--(10:0.5cm);          \draw (b2)--(-10:0.5cm);
                                          \end{scope}

                                                                                                                                             \begin{scope}[shift={(b21)}]
                    \draw (b21)--(80:0.5cm);          \draw (b21)--(100:0.5cm);
                                          \end{scope}

                                                                                                                                             \begin{scope}[shift={(b22)}]
                    \draw (b22)--(260:0.5cm);          \draw (b22)--(280:0.5cm);
                                          \end{scope}

                                                                                                                                             \begin{scope}[shift={(a3)}]
                    \draw (a3)--(260:0.5cm);          \draw (a3)--(280:0.5cm);
                                          \end{scope}

                                                                                                                                             \begin{scope}[shift={(a31)}]
                    \draw (a31)--(80:0.5cm);          \draw (a31)--(100:0.5cm);

                                                                                                                       \begin{scope}[shift={(b3)}]
                    \draw (b3)--(80:0.5cm);          \draw (b3)--(100:0.5cm);
                                          \end{scope}
                                          \end{scope}
    \node [below] at (2,-1) {$G$};

     \end{tikzpicture} \hskip2cm
     \begin{tikzpicture} [scale=0.5]
     \draw (2,2.3) ellipse (0.5cm and 3cm); \draw (6,2.3) ellipse (0.5cm and 3cm);
      \path    (2,0)  coordinate (e1);  \draw [fill=black] (e1) circle (0.07cm); \node [above] at (2,0) {\tiny $e_1$};   \node [above] at (6,0) {\tiny $e_2$};
            \path    (6,0)  coordinate (e2);  \draw [fill=black] (e2) circle (0.07cm); \draw (e1)--(e2);

            \path    (2,2)  coordinate (e3);  \draw [fill=black] (e3) circle (0.07cm);   \node [above] at (2,2) {\tiny $e_3$};
                   \path    (6,2)  coordinate (e5);  \draw [fill=black] (e5) circle (0.07cm);       \node [right] at (6,2) {\tiny $e_5$};
                               \path    (6,1.5)  coordinate (e4);  \draw [fill=black] (e4) circle (0.07cm);  \node [below] at (6,1.5) {\tiny $e_4$};
                        \path    (6,2.5)  coordinate (e6);  \draw [fill=black] (e6) circle (0.07cm);  \node [above] at (6,2.5) {\tiny $e_6$};
                        \path    (4,3)  coordinate (e10);  \draw [fill=black] (e10) circle (0.07cm);   \node [right] at (4,3) {\tiny $e_{10}$};
                        \draw (e3)--(e10) (e3)--(e4) (e3)--(e5) (e3)--(e6) (e10)--(e4) (e10)--(e5) (e10)--(e6) (e4)--(e6) (e4) .. controls (6.5,2) .. (e6);

                   \path    (2,4)  coordinate (e9);  \draw [fill=black] (e9) circle (0.07cm);    \node [below] at (2,4) {\tiny $e_9$};
                                        \path    (2,4.5)  coordinate (e7);  \draw [fill=black] (e7) circle (0.07cm);          \node [above] at (2,4.5) {\tiny $e_7$};
                                          \path    (6,4)  coordinate (e8);  \draw [fill=black] (e8) circle (0.07cm);              \node [above] at (6,4) {\tiny $e_8$};
                                          \draw (e10)--(e7) (e10)--(e8) (e10)--(e9) (e7)--(e8) (e7)--(e9) (e8)--(e9);

                                                                                                                           \begin{scope}[shift={(e1)}]
                    \draw (e1)--(260:0.5cm);          \draw (e1)--(280:0.5cm);
                                          \end{scope}

                                                                                                                                 \begin{scope}[shift={(e2)}]
                    \draw (e2)--(260:0.5cm);          \draw (e2)--(280:0.5cm);
                                          \end{scope}

                                                                                                                                         \begin{scope}[shift={(e3)}]
                    \draw (e3)--(170:0.5cm);          \draw (e3)--(190:0.5cm);
                                          \end{scope}

                                                                                                                                                                     \begin{scope}[shift={(e4)}]
                    \draw (e4)--(260:0.5cm);          \draw (e4)--(280:0.5cm);
                                          \end{scope}

                                                                                                                                                       \begin{scope}[shift={(e6)}]
                    \draw (e6)--(80:0.5cm);          \draw (e6)--(100:0.5cm);
                                          \end{scope}

                                                                                                                                                                     \begin{scope}[shift={(e5)}]
                    \draw (e5)--(10:0.5cm);          \draw (e5)--(-10:0.5cm);
                                          \end{scope}

                                                                                                                                                          \begin{scope}[shift={(e8)}]
                    \draw (e8)--(10:0.5cm);          \draw (e8)--(-10:0.5cm);
                                          \end{scope}

                                                                                                                                                                     \begin{scope}[shift={(e9)}]
                    \draw (e9)--(170:0.4cm);          \draw (e9)--(190:0.4cm);
                                          \end{scope}

                                                                                                                                                                     \begin{scope}[shift={(e7)}]
                    \draw (e7)--(170:0.35cm);          \draw (e7)--(190:0.35cm);
                                          \end{scope}

    \node [below] at (4,-1) {$L(G)$ };

    \draw [line width=1.2, blue] (e1)--(e2) (e3)--(e4)  (e3)--(e5) (e3)--(e6)  (e3)--(e10)  (e8)--(e7) (e8)--(e9) (e10)--(e7)  (e10)--(e9);

     \end{tikzpicture}

     Figure 3: Illustration of Definition \ref{BS}
\end{center}

By Observation \ref{Obs(G)}, we have the following.

\begin{Proposition}\label{edge cut}
Let $G\in {\cal G}(d,k)$ and
let $X$ be a minimal essential edge cut of $L(G)$ with
$V_{12}(X) = \{v_1, v_2, ..., v_m\}$. Each of the following holds.
\\
(i)If $V_1\neq \emptyset$ and $V_2\neq \emptyset$, then for each $j \in \{1,2\}$, the set $\bigcup\limits_{i=1}\limits^{m}A_{i,j,X}$ is an edge cut of $G$.
\\
(ii) For each $v_i\in V_{12}(X)$, $d_G(v_i)=a_{i,1,X}+a_{i,2,X}$.
\end{Proposition}

\section{Connectivity of second order iterated line graphs}

Recall that the graph family $\G(d,k)$ is defined in (\ref{def-g-k}). In this section, we need the following definition.
\begin{Definition}
Let $G\in {\cal G}(d,k)$. Fix a minimum essential edge cut $X(G)$  of $L(G)$, and let $m=|V_{12}(X(G))|$. We set $X=X(G)$ and
$V_{12}(X)=\{v_1,\cdots,v_m\}$. Using the notations in Definition \ref{BS}, we define
\begin{equation} \label{def-g1-4}
\begin{array}{ll}
\G_1=&\{G\in \G(d,k) \mbox{: $\exists i_0\in \{1, \ldots, m\}$ such
that $a_{i_0,1,X} = a_{i_0,2,X} = \lceil\frac{d}{2}\rceil$, $V_1(X) \neq \emptyset$, $V_2(X) \neq \emptyset$}\},
\\
\G_2=& \{G\in \G(d,k) \mbox{: $\forall i\in \{1,\ldots, m\}$
either $a_{i,1,X}\neq \lceil\frac{d}{2}\rceil$ or $a_{i,2,X} \neq \lceil\frac{d}{2}\rceil$, $V_1(X) \neq \emptyset$, $V_2(X) \neq \emptyset$} \},
\\
\G_3=&\{G\in \G(d,k) \mbox{:  exactly one of $V_1(X)$ and $V_2(X)$ is nonempty} \},
\end{array}
\end{equation}
\end{Definition}

By definition, each $\G_i$ depends on the values of $d$ and $k$, and $\G=\G(d,k)$ is disjoint union of $\G_1,\G_2$ and $\G_3$.
For each $i \in \{1,2,3\}$, define $ess'(L(\G_i)) =\min\{ ess'(L(G)) :G\in\G_i\}$.
Throughout this section, for each graph $G \in \G$, we fix one minimum essential edge cut $X(G)$ of $L(G)$. Thus for each graph $G \in \G$,
we use \eqref{V12} to define vertex subsets $V_1(X),V_2(X)$ and $V_{12}(X)$ of the graph $G$.

\begin{Lemma}\label{P1}
For $j\in\{1,2\}$, assume that $G\in {\cal G}_j$ is a graph with a minimum essential edge cut $X=X(G)$ of $L(G)$ satisfying $|X|=ess'(L(\G_j))$
such that $|E(G[V_{12}(X)])|$ is minimized,  among all minimum essential edge cuts of $L(G)$.
Denote $V_{12}(X)=\{v_1,v_2,\cdots,v_m\}$.
If $V_1(X)\not=\emptyset$ and $V_2(X)\neq \emptyset$, then each of the following holds.
\\
(i) $V_{12}(X)$  is a stable set of $G$.
\\
(ii) $\kappa'(G)=k$.
\\
(iii) If $d_G(v_{j_0})>d$ for some $v_{j_0}\in V_{12}(X)$, then $\sum\limits_{i=1}^m a_{i,1,X}=k$ and $a_{i,1,X}\le d$ for all $i$.
\\
(iv)  If $j=1$, then $d_G(v_i)=d$ for each $v_i\in V_{12}(X) - \{v_{i_0}\}$.
\\
(v) If $j=2$, then $d_G(v_i)=d$ for each $v_i\in V_{12}(X)$.
\end{Lemma}

\p  To prove each conclusions of Lemma \ref{P1}, we argue by contradiction to assume that
if one of these conclusions is false, then we are to find, for $j \in \{1,2\}$,
a graph $G' \in \G_j$ with $ess'(L(G')) < ess'(L(G))$, which leads to a contradiction of the assumption that
$ess'(L(G)) = ess'(L({\cal G}_j))$. Before the justifications of each conclusions, we assume that
if $G \in \G_j$ satisfies the hypotheses of this lemma, then we may assume that for the given
essential edge subset $X = X(G)$ in $L(G)$, the following is also satisfied:
\begin{equation} \label{local-completion}
\mbox{ for each $\ell \in \{1,2\}$, $|V_{\ell}(X)| \ge d+1$ and $G[V_{\ell}(X)]$ is a complete graph. }
\end{equation}
If fact, for each graph $G$ satisfying the hypotheses of the lemma, we obtain a new graph $\tilde{G}$ from $G$
by, for each $\ell \in \{1,2\}$, adding a set $W^{\ell}$ of new vertices with $|W^{\ell}| = \max\{0, d+1 - |V_{\ell}(X)|\}$ to $G$, and adding
a set $E^{\ell}$ of new edges such that $E^{\ell} \cup E(G[V_{\ell}(X)])$ is the edge set of a complete subgraph with vertex set $V_{\ell}(X) \cup W^{\ell}$
in $\tilde{G}$. Then by definition, if $G \in \G_j$, then $\tilde{G} \in \G_j$; and
the edge set $X$ in $L(G)$ is also a minimal essential edge cut of $L(\tilde{G})$,
with $E(G[V_{12}(X)]) = E(\tilde{G}[V_{12}(X)])$. Thus we can continue using $V_{12}(X)=\{v_1,v_2,\cdots,v_m\}$ in $\tilde{G}$.
By the choice of $W_{\ell}$, (\ref{local-completion}) holds.
Thus in the following, we always assume that $G$ satisfies (\ref{local-completion}).

To prove (i), we assume that  $V_{12}(X)$ is not a stable set. Then $E(G[V_{12}(X)]) \neq \emptyset$.
Construct a new graph $G'$ from $G$ by performing the following procedure (see Figure 4).

\noindent
\begin{tabular} {p{6.1in}}  \hline
{\bf Algorithm 1: } Construction of a new graph $G'$ from $G$ with given a given $X = X(G)$.  \\ \hline

1. {\bf Initialization.} Set $G: = G$, $E': = E(G[V_{12}(X)])$, $Q^1 = Q^2 = \emptyset$.
\\
2. {\bf Edge adjustment.}  If $E' = \emptyset$, then stop.
\\
If $E'$ contains an edge $e = v_iv_j$, then introduce  a new vertex $v_{ij} \notin V(G)$, update $V(G): = V(G) \cup \{v_{ij}\}$, $E(G): = E(G - e) \cup \{v_iv_{ij}, v_jv_{ij}\}$,
and $E':=E'-\{v_iv_j\}$.
\\
If $|E_{1,X}(v_i)|+|E_{1,X}(v_j)|\geq |E_{2,X}(v_i)|+|E_{2,X}(v_j)|$, then set
$Q^1: =  Q^1 \cup \{v_{ij}\}$, otherwise set $Q^2: =  Q^2 \cup \{v_{ij}\}$.
\\
3. {\bf Output.} A graph $G'$ in formed with $V(G'): = V(G) \cup Q^1 \cup Q^2$ and $E(G'): = E(G)$. \\ \hline
\end{tabular}

For each $e = v_iv_j \in E(G[V_{12}(X)])$, and for each
$f \in E(G)$ such that $ef \in X\subseteq E(L(G))$, we replace $ef$ in $X$ by $e'f$, where $e' = v_iv_{ij}$ if $f \in \partial_G(v_i)$ or
$e' = v_jv_{ij}$ if $f \in \partial_G(v_j)$. We call $e'f$ in $L(G')$ the corresponding edge of $ef$ in $L(G)$. For any edge $e_1e_2$
in $X$ such that neither $e_1$ nor $e_2$ is incident with a vertex in $\{v_i,v_j\}$, define the corresponding edge of $e_1e_2$ to be  itself.
Let $X'$ be formed from $X$ by changing each edge in $X$ with its corresponding edge. Then
$|X'| = |X|$.

By (\ref{local-completion}) and by Step 2 of Algorithm 1,
$X'$ is a minimal essential edge cut of $G'$ with $|V_{12}(X')| = |V_{12}(X)|$. Moreover, for $j \in \{1,2\}$, by Steps 2 and 3 of Algorithm 1 and by (\ref{local-completion}), if $G \in \G_j$,
then $G' \in \G_j$. As when Algorithm 1 stops, $E(G'[V_{12}(X)]) = \emptyset$, this contradicts the minimality of $|E(G[V_{12}(X)])|$. Thus (i) holds.

\begin{center}
\begin{tikzpicture} [scale=0.5]

\draw (0,2.3) ellipse (0.5cm and 3cm); \draw (4,2.3) ellipse (0.5cm and 3cm);   \node [above left] at (2.2,1.8) {\tiny $v_i$};   \node [above right] at (1.7,3.8) {\tiny $v_j$};

 \path    (2,0)  coordinate (v1);  \draw [fill=black] (v1) circle (0.07cm);  \path    (0,0)  coordinate (a1);  \draw [fill=black] (a1) circle (0.07cm);
  \path    (4,0)  coordinate (b1);  \draw [fill=black] (b1) circle (0.07cm); \draw (a1)--(b1);

  \path    (2,2)  coordinate (vi);  \draw [fill=black] (vi) circle (0.05cm);
    \path    (0,2)  coordinate (a2);  \draw [fill=black] (a2) circle (0.05cm); \draw (vi)--(a2);
      \path    (4,2)  coordinate (b2);  \draw [fill=black] (b2) circle (0.05cm); \draw (vi)--(b2);
        \path    (4,2.5)  coordinate (b21);  \draw [fill=black] (b21) circle (0.05cm); \draw (vi)--(b21);
        \path    (4,1.5)  coordinate (b22);  \draw [fill=black] (b22) circle (0.05cm); \draw (vi)--(b22);

   \path    (2,4)  coordinate (vj);  \draw [fill=black] (vj) circle (0.05cm);\draw (vi)--(vj);
      \path    (0,4)  coordinate (a3);  \draw [fill=black] (a3) circle (0.05cm);\draw (vj)--(a3);
      \path    (4,4)  coordinate (b3);  \draw [fill=black] (b3) circle (0.05cm);\draw (vj)--(b3);
      \path    (0,4.5)  coordinate (a31);  \draw [fill=black] (a31) circle (0.05cm);\draw (vj)--(a31);

    \node [above] at (1,-0.2) {\tiny $e_1$};      \node [above] at (3,-0.2) {\tiny $e_2$};

    \node [above] at (1,1.3) {\tiny $e_3$};          \node [above] at (3,1.3) {\tiny $e_4$};    \node [above] at (3.5,1.7) {\tiny $e_5$};    \node [above] at (3,2) {\tiny $e_6$};

    \node [above] at (1,3.5) {\tiny $e_9$};      \node [above] at (3,3.8) {\tiny $e_8$};     \node [above] at (1,4) {\tiny $e_7$};
      \node [right] at (1.8,3) {\tiny $e_{10}$};

                                                         \begin{scope}[shift={(a1)}]
                    \draw (a1)--(100:0.5cm);          \draw (a1)--(80:0.5cm);
                                          \end{scope}

                                                         \begin{scope}[shift={(b1)}]
                    \draw (b1)--(100:0.5cm);          \draw (b1)--(80:0.5cm);
                                          \end{scope}

                                                                                                   \begin{scope}[shift={(a2)}]
                    \draw (a2)--(170:0.5cm);          \draw (a2)--(190:0.5cm);
                                          \end{scope}

                                                                                                             \begin{scope}[shift={(b2)}]
                    \draw (b2)--(10:0.5cm);          \draw (b2)--(-10:0.5cm);
                                          \end{scope}

                                                                                                                                             \begin{scope}[shift={(b21)}]
                    \draw (b21)--(80:0.5cm);          \draw (b21)--(100:0.5cm);
                                          \end{scope}

                                                                                                                                             \begin{scope}[shift={(b22)}]
                    \draw (b22)--(260:0.5cm);          \draw (b22)--(280:0.5cm);
                                          \end{scope}

                                                                                                                                             \begin{scope}[shift={(a3)}]
                    \draw (a3)--(260:0.5cm);          \draw (a3)--(280:0.5cm);
                                          \end{scope}

                                                                                                                                             \begin{scope}[shift={(a31)}]
                    \draw (a31)--(80:0.5cm);          \draw (a31)--(100:0.5cm);

                                                                                                                       \begin{scope}[shift={(b3)}]
                    \draw (b3)--(80:0.5cm);          \draw (b3)--(100:0.5cm);
                                          \end{scope}
                                          \end{scope}
    \node [below] at (2,-1) {$G$};

     \end{tikzpicture} \hskip2cm
\begin{tikzpicture} [scale=0.5]

\draw (0,2.3) ellipse (0.5cm and 3cm); \draw (4,2.3) ellipse (0.5cm and 3cm);   \node [above left] at (2.2,1.8) {\tiny $v_i$};   \node [above right] at (1.7,3.8) {\tiny $v_j$};

 \path    (2,0)  coordinate (v1);  \draw [fill=black] (v1) circle (0.07cm);  \path    (0,0)  coordinate (a1);  \draw [fill=black] (a1) circle (0.07cm);
  \path    (4,0)  coordinate (b1);  \draw [fill=black] (b1) circle (0.07cm); \draw (a1)--(b1);

  \path    (2,2)  coordinate (vi);  \draw [fill=black] (vi) circle (0.05cm);
    \path    (0,2)  coordinate (a2);  \draw [fill=black] (a2) circle (0.05cm); \draw (vi)--(a2);
      \path    (4,2)  coordinate (b2);  \draw [fill=black] (b2) circle (0.05cm); \draw (vi)--(b2);
        \path    (4,2.5)  coordinate (b21);  \draw [fill=black] (b21) circle (0.05cm); \draw (vi)--(b21);
        \path    (4,1.5)  coordinate (b22);  \draw [fill=black] (b22) circle (0.05cm); \draw (vi)--(b22);

   \path    (2,4)  coordinate (vj);  \draw [fill=black] (vj) circle (0.05cm);  
      \path    (0,4)  coordinate (a3);  \draw [fill=black] (a3) circle (0.05cm);\draw (vj)--(a3);
      \path    (4,4)  coordinate (b3);  \draw [fill=black] (b3) circle (0.05cm);\draw (vj)--(b3);
      \path    (0,4.5)  coordinate (a31);  \draw [fill=black] (a31) circle (0.05cm);\draw (vj)--(a31);

   \path    (4,3.5)  coordinate (vij);  \draw [fill=black] (vij) circle (0.05cm); \draw (vij)--(vi) (vij)--(vj);   \node [right] at (3.8,3.5) {\tiny $v_{ij}$};

    \node [above] at (1,-0.2) {\tiny $e_1$};      \node [above] at (3,-0.2) {\tiny $e_2$};

    \node [above] at (1,1.3) {\tiny $e_3$};          \node [above] at (3,1.3) {\tiny $e_4$};    \node [above] at (3.5,1.7) {\tiny $e_5$};    \node [above] at (3,2) {\tiny $e_6$};

    \node [above] at (1,3.5) {\tiny $e_9$};      \node [above] at (3,3.8) {\tiny $e_8$};     \node [above] at (1,4) {\tiny $e_7$};

                                                         \begin{scope}[shift={(a1)}]
                    \draw (a1)--(100:0.5cm);          \draw (a1)--(80:0.5cm);
                                          \end{scope}

                                                         \begin{scope}[shift={(b1)}]
                    \draw (b1)--(100:0.5cm);          \draw (b1)--(80:0.5cm);
                                          \end{scope}

                                                                                                   \begin{scope}[shift={(a2)}]
                    \draw (a2)--(170:0.5cm);          \draw (a2)--(190:0.5cm);
                                          \end{scope}

                                                                                                             \begin{scope}[shift={(b2)}]
                    \draw (b2)--(10:0.5cm);          \draw (b2)--(-10:0.5cm);
                                          \end{scope}

                                                                                                                                             \begin{scope}[shift={(b21)}]
                    \draw (b21)--(80:0.5cm);          \draw (b21)--(100:0.5cm);
                                          \end{scope}

                                                                                                                                             \begin{scope}[shift={(b22)}]
                    \draw (b22)--(260:0.5cm);          \draw (b22)--(280:0.5cm);
                                          \end{scope}

                                                                                                                                             \begin{scope}[shift={(a3)}]
                    \draw (a3)--(260:0.5cm);          \draw (a3)--(280:0.5cm);
                                          \end{scope}

                                                                                                                                             \begin{scope}[shift={(a31)}]
                    \draw (a31)--(80:0.5cm);          \draw (a31)--(100:0.5cm);

                                                                                                                       \begin{scope}[shift={(b3)}]
                    \draw (b3)--(80:0.5cm);          \draw (b3)--(100:0.5cm);
                                          \end{scope}
                                          \end{scope}
    \node [below] at (2,-1) {$G'$};

     \end{tikzpicture}

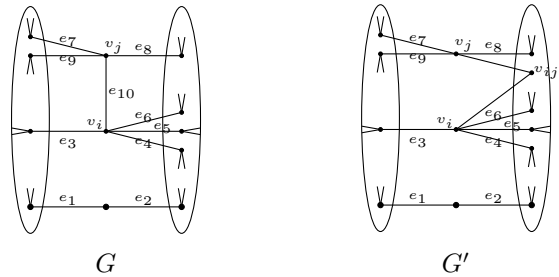
     Figure 4: Illustration of Step 2 for the essential edge cut $X=\{e_1e_2,e_3e_4,e_3e_5, e_3e_6,e_3e_{10}, e_8e_7, e_8e_9,$ $e_{10}e_7,e_{10}e_9\}$ in $L(G)$.
\end{center}

Since $G \in \G$, we have $\kappa'(G) \ge k$. By (i) and by Proposition \ref{edge cut},
$\sum\limits_{i=1}\limits^{m}a_{i,1, X} \ge k$, and if $\sum\limits_{i=1}\limits^{m}a_{i,1,X} = k$, then $\kappa'(G) = k$.
Since $\sum\limits_{i=1}\limits^{m}a_{i,1,X}\ge k$, for each $i$ with $1\leq i\leq m$, there exists an integer $a'_i$
with $0 < a'_i\leq a_{i,1,X}$ such that $\sum\limits_{i=1}\limits^{m} a'_i=k$. Construct a new graph $G'$ from $G$ by performing the following procedure.

\noindent
\begin{tabular}{p{6.1in}} \hline
{\bf Algorithm 2: } Construction of a new graph $G'$ from $G$ with a given $X = X(G)$.
\\ \hline

1. {\bf Initialization.}  Set $G := G$, and $T = Z^1 = Z^2 = \emptyset$.
Fix a sequence of positive integers $\{a_1', a_2', \ldots, a_m'\}$ satisfying $\sum\limits_{i=1}^m a_i' = k$
such that for each $i$, $t_i:= a_{i,1,X}-a_i' \ge 0$.
If $\max\{t_i: 1 \le i \le m\} = 0$, then stop.
\\
2.  {\bf Performing edge exchanges.} For each $i=1,2,\cdots,m$, if $t_i > 0$, then
choose an edge subset $T^i\subseteq  A_{i,1,X}$ with $|T^i|=t_i$ and a vertex subset $W^i \subseteq N_{i,2,X}(v_i)$ with $|W^i|=t_i$. Update $T: = T \cup T^i$.
\\
If $|E_{1,X}(v_i)|>|E_{2,X}(v_i)|$, then set $Z^1:=Z^1 \cup  \{v_iy: y \in W^i\}$; if $|E_{1,X}(v_i)|\le |E_{2,X}(v_i)|$, then set $Z^2:=Z^2 \cup  \{v_iy: y \in W^i\}$.
\\
3. {\bf Output.} A graph $G'$ with $V(G'): = V(G)$ and $E(G'): = (E(G) \cup Z^1 \cup Z^2) - T$.
\\ \hline
\end{tabular}

In the $i$th iteration within Step 2 of Algorithm 2, when $t_i > 0$, we denote $T_i = \{f_1^i, f_2^i, ..., f_{t_i}^i\}$ and
$W_i = \{w_1^i, w_2^i, ..., w_{t_i}^i\}$. Fix an $i$ with $1 \le i \le m$. If for some edge $e \in E(G)$ and some $f_{\ell}^i \in T_i$,
an edge $ef_{\ell}^i \in X \subseteq E(L(G))$,
then we define $f' = v_iw_{\ell}^i$, and the corresponding edge of $ef_{\ell}$ to be $ef_{\ell}'$.
If an edge $ef$ of $E(L(G))$ is not of the form $ef_{\ell}^i$ as described above, then define
the corresponding edge of $ef$ to be itself.
Let $X'$ be the set obtained from $X$ by changing all edges in $X$ to their corresponding edges.
Then by (\ref{local-completion}) and by Steps 2 and 3 of Algorithm 2, $X'$ is a minimal an essential edge cut of $L(G')$ with $V_{12}(X)=V_{12}(X')$.
By (i), $V_{12}(X)$ is a stable set in $G$.
As in Step 3, all added edges will not joint two vertices in $V_{12}(X')$, we note that
$V_{12}(X')$ is a stable set in $G'$.
For $j=1,2$, Algorithm 2 and (\ref{local-completion}) also ensure  that  if $G\in {\cal G}_j$, then $G'\in {\cal G}_j$.
By Definition \ref{BS} and Step 2 of Algorithm 2, $a_{i,2,X}\ge a_{i,1,X} \ge a_{i,1,X'}$.
By Proposition \ref{edge cut}, $d_G(v_i)=a_{i,1,X}+a_{i,2,X} \ge a_{i,1,X}+a_{i,1,X'}$.
As $|X|=ess'(L({\cal G}_j))$, by \eqref{eq} we are lead to a contradiction:
\begin{eqnarray*}
0 \ge |X|-|X'|
&= & \sum\limits_{i=1}\limits^{m}a_{i,1,X} \cdot(d_{G}(v_i)-a_{i,1,X})-
\sum\limits_{i=1}\limits^{m} a_{i,1,X'} \cdot (d_{G'}(v_i)-a_{i,1,X'})\\
&= & \sum\limits_{i=1}\limits^{m} \Big[d_G(v_i) (a_{i,1,X}-a_{i,1,X'})-(a_{i,1,X}^2-a_{i,1,X'}^2)\Big] \\
&=& \sum\limits_{i=1}^m (a_{i,1,X}-a_{i,1,X'}) (d_G(v_i)-a_{i,1,X}-a_{i,1,X'})>0.
\end{eqnarray*}
This proves  (ii).

To prove (iii), we by contradiction assume that $\sum\limits_{i=1}^m a_{i,1,X}\ge k+1$.
By Proposition \ref{edge cut}, $a_{j_0,1, X}+a_{j_0,2,X} = d_G(v_{j_0})>d \ge 3$.
Thus, either $a_{j_0,1,X}\ge 2$ or both $a_{j_0,1,X}=1$ and $a_{j_0,2,X}\ge d\ge 3$.

If $a_{j_0,1,X}\ge 2$, then pick an edge $e' \in A_{j_0,1,X}$, and
if $a_{j_0,1,X} = 1$, then pick an edge $e' \in A_{j_0,2,X}$. Once $e'$ is chosen, let $G' = G-\{e'\}$.
We consider how the deletion of $e'$ will change the essential edge cut $X$ of $L(G)$ to an edge cut $X'$ of $L(G')$.
By the formation of $G'$, let $X' = X-Y^1$, where
\[
Y^1 =
\left\{
\begin{array}{ll}
\{e'f\in X: f \in A_{j_0,1,X}\} & \mbox{ if $a_{j_0,1,X}\geq 2$,}
\\
\{e'f\in X: f \in A_{j_0,2,X}\} & \mbox{ if $a_{j_0,1, X}=1$.}
\end{array} \right.
\]
As $X$ is a minimal essential edge cut of $L(G)$ and (\ref{local-completion}) holds, $X'$ is a minimal essential edge cut of $L(G')$ with $V_{12}(X')=V_{12}(X)$.
Note that
$$
(a_{j_0,1,X'}, a_{j_0,2, X'})=\begin{cases}
  (a_{j_0,1, X'}-1, a_{j_0,2,X'})    & \text{ if } a_{j_0,1,X}\ge 2, \\
   (a_{j_0,1, X'}, a_{j_0,2, X'}-1)   & \text{ if } a_{j_0,1, X}=1,
\end{cases}
$$
and so $d_{G'}(v_{j_0})=d_G(v_{j_0})-1\ge d$, and $\sum\limits_{i=1}^m a_{i,1,X'}\ge \sum\limits_{i=1}^m a_{i,1,X}-1\ge (k+1)-1=k$.
By Proposition \ref{edge cut}  and by (\ref{local-completion}), for $j \in \{1,2\}$, if $G \in \G_j$, then $G' \in \G_j$. As $|X|=ess'(L({\cal G}_j))$, by \eqref{eq},
the following contradiction is obtained:
$$
0 \ge |X|-|X'|= a_{j_0,1,X}a_{j_0,2, X}-a_{j_0,1, X'}a_{j_0,2,X'} =\begin{cases}
    a_{j_0,2,X},  & \text{ if } a_{j_0,1, X}\ge 2 \\
    a_{j_0,1, X},  & \text{ if } a_{j_0,1, X}=1
\end{cases}
\; \; > 0.
$$
Therefore, we must have $\sum\limits_{i=1}^m a_{i,1,X}=k\le d$.
This implies that $a_{i,1,X}\le d$ for all $i$, and so (iii) holds.

In the following, we prove (iv) and (v). Pick a vertex $v_i\in V_{12}(X)$. Then $d_G(v_i)=a_{i,1,X}+a_{i,2,X}\ge d$ and $a_{i,1,X}\le a_{i,2,X}$.
Thus $a_{i,2,X} \ge \frac{a_{i,1,X}+a_{i,2,X}}{2} \ge \lceil\frac{d}{2}\rceil$.
First we claim that \begin{equation} \label{3-1}
\hbox{ if  there exists }
i_0\in\{1,2,\cdots,m\} \hbox{  such that } a_{i_0,1,X}=\Big\lceil\frac{d}{2}\Big\rceil,
\hbox{then } a_{i_0,2,X}=\Big\lceil\frac{d}{2}\Big\rceil.
\end{equation}
Define $t_{i_0} = a_{i_0,2,X}- \lceil\frac{d}{2}\rceil$. Assume that (\ref{3-1}) is false. Then $t_{i_0} > 0$.
Choose a subset $T^{i_0}\subseteq  A_{i_0,2,X}$ with $|T^{i_0}|=t_{i_0}$ and define $G' = G - T^{i_0}$ and
$X' = X- \bigcup_{e \in T^{i_0}} \{ef \in E(L(G)): f \in E(G)\}$.
As $X$ is a minimal essential edge cut of $L(G)$ and by (\ref{local-completion}), $X'$ is a minimal essential edge cut of $L(G')$
with $V_{12}(X')=V_{12}(X)$ and with $\sum\limits_{i=1}^{m} |A_{i,1,X'}|=\sum\limits_{i=1}^m |A_{i,1,X}|\ge k$.
By the formation of $G'$, we have, $a_{i,1,X}=a_{i,1,X'}$ and $a_{i,2,X}=a_{i,2,X'}$ for each $i\not=i_0$,
while $a_{i_0,1,X}=a_{i_0,1,X'}=\lceil \frac{d}{2}\rceil$ and $a_{i_0,2,X'}=\lceil \frac{d}{2}\rceil$.
By (\ref{local-completion}), it follows that for $j=1,2$, if $G\in {\cal G}_j$, then $G'\in {\cal G}_j$. As $|X|=ess'(L({\cal G}_j))$ and $t_{i_0} > 0$,
we have a contradiction:
$0 \ge |X|-|X'| \ge |X| - (|X| - t_{i_0}) > 0$,
which proves (\ref{3-1}).

Next we claim that
\begin{equation} \label{3-2}
\hbox{for each } i\in\{1,2,\cdots,m\},
\hbox{if }
d_G(v_i)>d,
\hbox{then } a_{i,1,X}\ge \Big\lceil\frac{d}{2}\Big\rceil.
\end{equation}
Let $t_{i}'=d_G(v_{i})-d$, where $1 \le i \le m$.
Suppose that (\ref{3-2}) is false. Then there exists an index $i_0\in\{1,2,\cdots,m\}$ with both $t_{i_0}'>d$
and $a_{i_0,1,X}\le \lfloor\frac{d}{2}\rfloor$.
Choose a subset $T^{i_0'}\subseteq  A_{i_0,2,X}$ with $|T^{i_0'}|=t_{i_0}'$ and define $G' = G - T^{i_0'}$ and
$X' = X- \bigcup_{e \in T^{i_0'}} \{ef \in E(L(G)): f \in E(G)\}$.
As $X$ is a minimal essential edge cut of $L(G)$ and by (\ref{local-completion}), $X'$ is a minimal essential edge cut of $L(G')$
with $V_{12}(X')=V_{12}(X)$.
As $a_{i_0,2,X}-(d_G(v_{i_0})-d)=d-a_{i_0,1,X}\ge \lceil\frac{d}{2}\rceil\ge a_{i_0,1,X}$, we have  $a_{i,1,X}=a_{i,1,X'}$ and $a_{i,2,X}=a_{i,2,X'}$ for each $i\not=i_0$, while $a_{i_0,1,X}=a_{i_0,X'}$ and $a_{i_0,2,X'}=d-a_{i_0,1,X}$. It follows that $\sum\limits_{i=1}^{m} |A_{i,1,X'}|=\sum\limits_{i=1}^m |A_{i,1,X}|\ge k$.
Hence  by (\ref{local-completion}),
for $j=1,2$, if $G\in {\cal G}_j$, then $G'\in {\cal G}_j$. As $|X|=ess'(L({\cal G}_j))$ and $t_{i_0} > 0$,
we have a contradiction:
$0 \ge |X|-|X'| \ge |X| - (|X| - t_{i_0}') > 0$,
which proves (\ref{3-2}).

In order to prove (iv) and (v), we consider two cases depending on whether  there exists $i_0\in\{1,2,\cdots,m\}$ such that $a_{i_0,1,X}=a_{i_0,2,X}=\lceil\frac{d}{2}\rceil$ or not.
Define $$S=\left\{\begin{array}{ll}
\{1,2,\cdots,m\}-\{i_0\}, & \hbox{ if such an } i_0 \hbox{ exists,} \\
\{1,2,\cdots,m\}, & \hbox{ otherwise.}
\end{array}\right.$$
We are to prove that $d_G(v_i)=d$ for every $i\in S$  in order to complete the proofs of  (iv) and (v).
Suppose, to the contrary, that there exists $j_0\in S$ such that $d_G(v_{j_0})> d$. Let $s= a_{j_0,1,X}-\lceil \frac{d}{2}\rceil$.
We are to show that $s > 0$.

By (\ref{3-2}), $a_{j_0,1,X}\ge \lceil\frac{d}{2}\rceil$. By (iii), $\sum\limits_{i=1}^m a_{i,1,X}=k$.
If $s = 0$, then $a_{j_0,1, X}= \lceil\frac{d}{2} \rceil$, and so by (\ref{3-1}), $a_{j_0,2,X}=\lceil\frac{d}{2} \rceil$.
Hence by Proposition \ref{edge cut}, $2  \lceil\frac{d}{2} \rceil = a_{j_0,1, X} + a_{j_0,2,X} = d_G(v_{j_0})>d$,
forcing that $d$ is odd and $a_{j_0,1,X}+a_{j_0,2,X}= d+1$. By the definition of $S$, $j_0\not=i_0$.
As $a_{i_0,1,X}=a_{j_0,2, X}=\lceil\frac{d}{2}\rceil$, a contradiction is obtained:
$$d+1=a_{j_0,1,X}+a_{j_0,2,X}=a_{j_0,1,X}+a_{i_0,1,X} \le \sum\limits_{i=1}^m a_{i,1,X}=k\le  d. $$

Therefore, $s > 0$, and so $a_{j_0,2,X} \ge a_{j_0,1,X}> \lceil\frac{d}{2}\rceil$.
Choose edge subsets $T_1^{j_0}\subseteq A_{j_0,1, X}$ and $T_2^{j_0}\subseteq A_{j_0,2,X}$ with $|T_1^{j_0}|=|T_2^{j_0}|=s$.
By (\ref{local-completion}), for $\ell \in \{1,2\}$, each $|V_\ell(X)| \ge d+1$. Pick vertex subsets $Q_1^{m+1} \subseteq V_1(X)$ and $Q_2^{m+1} \subseteq V_2(X)$
with $|Q_1^{m+1}|=s$ and $|Q_2^{m+1}|=d-s$. Let $v_{m+1}$ be a new vertex not in $V(G)$.
Obtain a new graph $G'$ with
\begin{equation} \label{alg-6}
\begin{array}{ll}
V(G') = & V(G) \cup \{ v_{m+1}\} \mbox{ and }
\\
E(G')  = & \left(E(G) \cup  \{v_{m+1}z: z\in Q_1^{m+1}\} \cup \{v_{m+1}z: z\in Q_2^{m+1}\} \right) - (T_1^{j_0} \cup T_2^{j_0}).
\end{array}
\end{equation}
Let $s'$ be an integer with $a_{j_0,2,X}=\lceil\frac{d}{2}\rceil+s'$. Since $a_{j_0,1,X}=\Big\lceil\frac{d}{2}\Big\rceil+s$,  we have $s' \ge s$.
By (\ref{alg-6}), $d_{G'}(v_{j_0})=\lceil\frac{d}{2}\rceil +\big(\lceil\frac{d}{2}\rceil+s'-s\big)\ge d$
and $d_{G'}(v_{m+1})=s+(d-s)=d$. Define
\[
X' = X \cup  \{f_1f_2: f_i=v_{m+1}z_i,  z_i \in Q_i^{m+1}, 1 \le i \le 2\} - \{e'f': e' \in T_1^{j_0} \cup T_2^{j_0} \mbox{ and } f' \in E(G)\}.
\]
By (\ref{local-completion}) and (\ref{alg-6}), $X'$ is a minimal essential edge cut of $L(G')$, with
$a_{j_0,1, X'}=\lceil\frac{d}{2}\rceil$ and $a_{j_0,2,X'}=\lceil\frac{d}{2}\rceil+s'-s$.
Since $d_G(v_{j_0})=a_{j_0,1,X}+a_{j_0,2,X}\ge d$,  it follows by (iii) that  $\lceil\frac{d}{2}\rceil<a_{j_0,1,X}\le d$. Hence,
$s=a_{j_0,1,X}-\lceil\frac{d}{2}\rceil\le \lfloor \frac{d}{2}\rfloor$, and so $s\le d-s$. This implies that $a_{m+1,1, X'}=s$ and $a_{m+1,2,X'}=d-s$.
Since $a_{j_0,1, X}=\lceil\frac{d}{2}\rceil+s$, it follows that
$\sum\limits_{i=1}^{m+1} a_{i,1,X'}=\sum\limits_{i=1}^{m} a_{i,1, X}\ge k$. Thus, for $j=1,2$, if $G\in {\cal G}_j$, then $G'\in {\cal G}_j$.
As $|X|=ess'(L({\cal G}_j))$, by \eqref{eq},
\begin{eqnarray*}
0 \ge |X|-|X'| &=& a_{j_0,1,X}a_{j_0,2,X} -[a_{j_0,1,X'}a_{j_0,2,X'} +a_{m+1,1,X'} a_{m+1,2,X'} ] \\
&=& \Big(\Big\lceil\frac{d}{2}\Big\rceil+s\Big)\Big(\Big\lceil\frac{d}{2}\Big\rceil+s'\Big)
-\Big[\Big\lceil\frac{d}{2}\Big\rceil \cdot \Big(\Big\lceil\frac{d}{2}\Big\rceil+s'-s\Big)+
s(d-s)\Big] \\
&=&  2s \Big\lceil\frac{d}{2}\Big\rceil+ss'-sd+s^2 \ge ss'+s^2>0.
\end{eqnarray*}
This contradiction indicates that for any $i \in S$, we always have $d_G(v_{i}) = d$, and so (iv) and (v) hold.
\q

Now, we consider graphs in $\G_3$. For each graph $G \in \G_3$ with $X = X(G)$, by (\ref{def-g1-4}), exactly one of $V_1(X)$ and $V_2(X)$ is nonempty. By symmetry, in the discussions below, we always assume that $V_2(X)=\emptyset$ in this case,  and so
for each $i$ with $1 \le i \le m = |V_{12}(X)|$, both $A_{i,1,X}=\emptyset$ and $A_{i,2,X}=E_{G}[v_i,V_1(X)]$.


\begin{Lemma}\label{P3}
Let $G\in \G_3\cap \G_0(d,k)$ with a minimum edge cut $X=X(G)$ with $|X|=ess'(L(\G_3))$, and with
$V_1(X)\not=\emptyset$ and $V_2(X)=\emptyset$,
such that $|E(G[V_{12}(X)])|$ is minimized among all minimum essential edge cuts of $L(G)$.
Denote $V_{12}(X)=\{v_1,v_2,\cdots,v_m\}$ and let $f_X$  denote the $X$-induced 2-edge-coloring of $G$ as in Definition \ref{v-x}.
Each of the following holds.
\\
(i) Let $E'=\{v_iv_j\in E(G[V_{12}(X)]: f_X(v_{i}v_{j})=1\}$. Then $E'=\emptyset.$
\\
(ii) $d_G(v_i)=d$ for $1\leq i\leq m$.\\
(iii)  For some integer $m$ with $2\leq m\leq k$,  $G[V_{12}(X)] \cong K_{1,m-1}$.
\end{Lemma}

\p Suppose that (i) is false and let $e_0 = v_{i_0}v_{j_0} \in E'$. Pick a vertex subset $W^1$ disjoint from $V(G)$ with
$|W^1| = \max\{0, d+1-|V_1(X)|\}+1$, a distinguished vertex $v_{i_0j_0}  \in W^1$ and define a new graph $G'$ (see Figure 5) as follows:
\begin{equation} \label{alg-7}
\begin{array}{ll}
V(G') = & V(G) \cup W^1,
\\
E(G') = & (E(G) - \{v_{i_0}v_{j_0} \}) \cup \{v_{i_0}v_{i_0j_0}, v_{i_0j_0}v_{j_0} \} \cup \{vw: v \in V_1(X) \mbox{ and } w \in W^1\}.
\end{array}
\end{equation}
Thus $G'[W^1 \cup V_1(X)]$ is a complete subgraph of $G'$ with order at least $d+1$.
As $G \in \G(d,k)$, and by (\ref{alg-7}), we also have $G' \in \G(d,k)$.
For each edge $ee' \in X$, if $e=e_0$ and $e' \in E_{2,X}(v_{i_0})\cup E_{2,X}(v_{j_0})$, then
define the corresponding edge of $e_0e'$ to be $v_{i_0}v_{i_0j_0}e'$ when $V(e_0) \cap V(e') = \{v_{i_0}\}$
or $v_{j_0}v_{i_0j_0}e'$ when $V(e_0) \cap V(e') = \{v_{j_0}\}$; otherwise, define the corresponding edge of $ee'$
to be itself.
 Let $X'$ be formed from $X$ by changing each edge in $X$ to its corresponding edge.
As $X$ is a minimal essential edge cut of $L(G)$ and by (\ref{alg-7}), $X'$ is a minimal essential edge cut of $L(G')$ with $V_{12}(X)=V_{12}(X')$,
$V_2(X')=\emptyset$, and $|X'|=|X|$. Thus, $G' \in  \G_3\cap {\cal G}_0(d,k)$. Since $|E(G[V_{12}(X)])| = |E(G'[V_{12}(X')])|+1$,
a contradiction to the minimality of $|E(G[V_{12}(X)])|$ is obtained, which implies (i).

\begin{center}
\begin{tikzpicture} [scale=0.5]

\draw (0,2.3) ellipse (0.6cm and 3cm);
\node [right] at (1.8,2) {\tiny $v_2$};   \node [above] at (2,3.8) {\tiny $v_3$};
\node [below] at (2,0) {\tiny $v_1$};

 \path    (2,0)  coordinate (v1);  \draw [fill=black] (v1) circle (0.07cm);  \path    (0,0)  coordinate (a1);  \draw [fill=black] (a1) circle (0.07cm);
  \path    (2,2)  coordinate (vi);  \draw [fill=black] (vi) circle (0.05cm);
    \path    (0,2)  coordinate (a21);  \draw [fill=black] (a21) circle (0.05cm); \draw (vi)--(a21);
   \path    (2,4)  coordinate (vj);  \draw [fill=black] (vj) circle (0.05cm);\draw (vi)--(vj);
         \path    (0,4)  coordinate (a3);  \draw [fill=black] (a3) circle (0.05cm);\draw (vj)--(a3);

\draw (a1)--(v1); \draw (a3)--(vj); \draw (v1)--(vj);

    \node [above] at (1,-0.2) {\tiny $e_1$};

                 \node [above] at (1,1.9) {\tiny $e_2$}; \node [above] at (1,3.9) {\tiny $e_3$};
                       \node [left] at (2.1,3) {\tiny $e_{6}$};             \node [left] at (2.1,1.1) {\tiny $e_{5}$};

\draw (vj) .. controls (3,2).. (v1);  \node [right] at (2.5,2.3) {\tiny $e_4$};

  \begin{scope}[shift={(a1)}]
\draw (a1)--(170:0.5cm);          \draw (a1)--(190:0.5cm);
\end{scope}

  \begin{scope}[shift={(a21)}]
\draw (a21)--(170:0.5cm);          \draw (a21)--(190:0.5cm);
 \end{scope}

 \begin{scope}[shift={(a3)}]
\draw (a3)--(170:0.5cm);          \draw (a3)--(190:0.5cm);
  \end{scope}

    \node [below] at (1,-1) {\tiny $V_{12}(X)=\{v_1,v_2,v_3\}$};

     \end{tikzpicture} \hskip1cm
    \begin{tikzpicture} [scale=0.5]

\draw (0,2.3) ellipse (0.6cm and 3cm);

    \path    (0,0)  coordinate (e1);  \draw [fill=black] (e1) circle (0.05cm);  \node [above] at (e1) {\tiny $e_1$};
         \path    (0,2)  coordinate (e2);  \draw [fill=black] (e2) circle (0.05cm);\node [above] at (e2) {\tiny $e_2$};
     \path    (0,3.9)  coordinate (e3);  \draw [fill=black] (e3) circle (0.05cm);\node [above] at (e3) {\tiny $e_3$};

        \path    (3,0)  coordinate (e4);  \draw [fill=black] (e4) circle (0.05cm);  \node [right] at (e4) {\tiny $e_4$};
         \path    (3,2)  coordinate (e5);  \draw [fill=black] (e5) circle (0.05cm); \node [right] at (e5) {\tiny $e_5$};
     \path    (3,3.9)  coordinate (e6);  \draw [fill=black] (e6) circle (0.05cm); \node [right] at (e6) {\tiny $e_6$};


\draw (e1)--(e4) (e1)--(e5) (e2)--(e5) (e2)--(e6) (e3)--(e4) (e3)--(e6);
\draw (e4)--(e6) (e4) .. controls (4,2)..(e6);

                                                         \begin{scope}[shift={(e1)}]
                    \draw (e1)--(170:0.5cm);          \draw (e1)--(190:0.5cm);
                                          \end{scope}

                                                                                                   \begin{scope}[shift={(e2)}]
                    \draw (e2)--(170:0.5cm);          \draw (e2)--(190:0.5cm);
                                          \end{scope}

                                                                                                   \begin{scope}[shift={(e3)}]
                    \draw (e3)--(170:0.5cm);          \draw (e3)--(190:0.5cm);
                                          \end{scope}

\draw [blue, line width=1.2] (e5)--(e6) (e4) .. controls (4,2)..(e6) (e2)--(e5) (e1)--(e5) (e1)--(e4) (e3)--(e4);

    \node [below] at (1,-1) {\tiny $L(G)$};

     \end{tikzpicture} \hskip2cm
     \begin{tikzpicture} [scale=0.5]

\draw (0,2.3) ellipse (0.6cm and 3cm);
\node [right] at (1.8,2) {\tiny $v_2$};   \node [above] at (2,3.8) {\tiny $v_3$};
\node [below] at (2,0) {\tiny $v_1$};

 \path    (0,0)  coordinate (a1);  \draw [fill=black] (a1) circle (0.07cm);
     \path    (0,1.5)  coordinate (a21);  \draw [fill=black] (a21) circle (0.05cm);
     \path    (0,2.5)  coordinate (a22);  \draw [fill=black] (a22) circle (0.05cm);  \node [left] at (0.2,2.5) {\tiny $v_{23}$};
                   \path    (0,4)  coordinate (a3);  \draw [fill=black] (a3) circle (0.05cm);

 \path    (2,0)  coordinate (v1);  \draw [fill=black] (v1) circle (0.07cm);
  \path    (2,2)  coordinate (vi);  \draw [fill=black] (vi) circle (0.05cm);
   \path    (2,4)  coordinate (vj);  \draw [fill=black] (vj) circle (0.05cm);
\draw (a22)--(vj) (a22)--(vi); \draw (a21)--(vi);
\draw (vj)--(a3);
\draw (a1)--(v1); \draw (a3)--(vj); \draw (v1)--(vi);
\draw [dotted] (vi)--(vj);

    \node [above] at (1,-0.2) {\tiny $e_1$};

                 \node [above] at (1,1.2) {\tiny $e_2$}; \node [above] at (1,3.9) {\tiny $e_3$};
                           \node [left] at (2.1,1.1) {\tiny $e_{5}$};

\draw (vj) .. controls (3,2).. (v1);  \node [right] at (2.5,2.3) {\tiny $e_4$};

  \begin{scope}[shift={(a1)}]
\draw (a1)--(170:0.5cm);          \draw (a1)--(190:0.5cm);
\end{scope}

  \begin{scope}[shift={(a21)}]
\draw (a21)--(170:0.5cm);          \draw (a21)--(190:0.5cm);
 \end{scope}

 \begin{scope}[shift={(a3)}]
\draw (a3)--(170:0.5cm);          \draw (a3)--(190:0.5cm);
  \end{scope}

    \node [below] at (1,-1) {\tiny $G'$};

     \end{tikzpicture} \hskip1cm
     \begin{tikzpicture} [scale=0.5]

\draw (0,2.3) ellipse (0.6cm and 3cm);

    \path    (0,0)  coordinate (e1);  \draw [fill=black] (e1) circle (0.05cm);  \node [above] at (e1) {\tiny $e_1$};
         \path    (0,1)  coordinate (e2);  \draw [fill=black] (e2) circle (0.05cm);\node [above] at (e2) {\tiny $e_2$};
     \path    (0,2)  coordinate (e-1);  \draw [fill=black] (e-1) circle (0.05cm);\node [left] at (e-1) {\tiny $v_{23}v_2$};
\path    (0, 3)  coordinate (e-2);  \draw [fill=black] (e-2) circle (0.05cm);  \node [left] at (e-2) {\tiny $v_{23}v_3$};
\path    (0,4)  coordinate (e3);  \draw [fill=black] (e3) circle (0.05cm);  \node [right] at (e3) {\tiny $e_3$};

         \path    (3,1)  coordinate (e5);  \draw [fill=black] (e5) circle (0.05cm); \node [right] at (e5) {\tiny $e_5$};
     \path    (3,3)  coordinate (e4);  \draw [fill=black] (e4) circle (0.05cm); \node [right] at (e4) {\tiny $e_4$};

\draw (e1)--(e4) (e1)--(e5) (e2)--(e5) (e-1)--(e5) (e-2)--(e3) (e-2)--(e4) (e3)--(e4) (e4)--(e5);

                                                         \begin{scope}[shift={(e1)}]
                    \draw (e1)--(170:0.5cm);          \draw (e1)--(190:0.5cm);
                                          \end{scope}

                                                                                                   \begin{scope}[shift={(e2)}]
                    \draw (e2)--(170:0.5cm);          \draw (e2)--(190:0.5cm);
                                          \end{scope}

                                                                                                   \begin{scope}[shift={(e3)}]
                    \draw (e3)--(170:0.5cm);          \draw (e3)--(190:0.5cm);
                                          \end{scope}


    \node [below] at (1,-1) {\tiny $L(G')$};

     \end{tikzpicture}

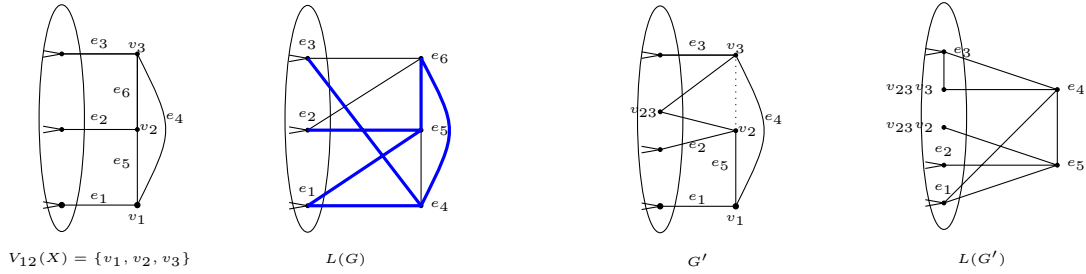
Figure 5: Illustration of forming $G'$ in (\ref{alg-7}) with $X=\{e_1e_4,e_1e_5,e_2e_5,e_3e_4,e_5e_6,e_4e_6\}$, $f_X(e_6)=1$ and $i_0 = 2$, $j_0 = 3$.
\end{center}

We argue by contradiction to prove (ii) and assume that $d_{G}(v_{i_0})>d$ for some $i_0$. Let $t_{i_0} = d_G(v_{i_0})-d$.
As $G \in \G_0(d,k)$ and by Proposition \ref{m}, $m\le k$.
By Observation \ref{Obs(G)}(iii), $d_G(v_{i_0})=|E_{1,X}(v_{i_0})|+|E_{2,X}(v_{i_0})|$.
If $|E_{1,X}(v_{i_0})|\le d_G(v_{i_0})-d$, then $|E_{2,X}(v_{i_0})|\ge  d$, leading to a contradiction $k \ge m\ge |E_{2,X}(v_{i_0})| +1\ge d+1\ge k+1$.
Thus $|E_{1,X}(v_{i_0})| > d_G(v_{i_0})-d$. By (i), $\partial_G(V_{12}(X)) \supseteq  E_{1,X}(v_{i_0})=\emptyset$.
We are to construct a new graph $G' \in \G_3 \cap \G_0(d,k)$ to find a contradiction.

Choose a set $T^{i_0}\subseteq E_{1,X}(v_{i_0})$ with $|T^{i_0}|=t_{i_0} > 0$, and a set
$W^1$ of vertices disjoint from $V(G)$ with $|W^1| = \max\{0, d+1-|V_1(X)|\}$. Define a graph $G'$ with
\begin{equation} \label{alg-8}
\begin{array}{ll}
V(G') = & V(G) \cup W^1,
\\
E(G') = & (E(G) - T^{i_0}) \cup\{vw: v \in V_1(X) \mbox{ and } w \in W^1\}.
\end{array}
\end{equation}
As $V_2(X)=\emptyset$, $V(G')$ is a disjoint union of $V_{12}(X)$ and $W^1 \cup  V_1(X)$.
To show that $\kappa'(G')\geq k$,
we randomly pick a proper nonempty subset $V' \subset V(G')$. If both $V'\cap (W^1 \cup  V_1(X)) \neq\emptyset$ and $V' - (W^1 \cup  V_1(X)) \neq \emptyset$,
then by (\ref{alg-8}),  $G'[W^1 \cup  V_1(X)]$ is a complete graph of order at least $d+1 \ge k+1$, which implies that
$|\partial_{G'}(V')| \geq \kappa'(G'[W^1 \cup  V_1(X)]) \ge d \ge k$.
Hence we assume that for some $t$ with $1 \le t \le m$,  $V'=\{v_1,\ldots,v_t\} \subseteq V_{12}(X)$.
Define
$$
f(t)=|\partial_G(\{v_1,\ldots,v_t\})|-k, \mbox{ and } q(t) = dt-\frac{t(t-1)}{2} -k =-\frac{t^2}{2}+\frac{2d+1}{2}t-k.
$$
By definition, $f(t) \ge q(t)$. Direct computation yields that
$q(1)\ge d-k\geq 0$ and $q(k)=-\frac{k^2}{2}+\frac{2d+1}{2}k-k=\frac{k}{2}(2d-1-k)\geq 0$.
As $q(t)$ is a quadratic function in $t$ on $[1,k]$ with a negative leading coefficient, and $d\ge k\ge 1$,
we have $q(t) \ge 0$ on $[1,k]$. This implies that $|\partial_G(V')|-k = f(t) \ge q(t) \ge 0$, for any $t$ with $1 \le t\le m<k$.
This, by definition, implies that $\kappa'(G')\geq k$.

Define $Y^3 = \{ef \in X: e\in T^{i_0}, f\in E_{2,X}(v_{i_0})\}$, and let $X' = X - Y^3$.
By (\ref{alg-8}), $X'$ is a minimal essential edge cut of $L(G')$ with $V_{12}(X)=V_{12}(X')$.
Again by (\ref{alg-8}),  $d_{G'}(v_{i_0})=d$,  $V_2(X')=V_2(X)=\emptyset$ and $V_1(X') \supseteq V_1(X) \not=\emptyset$.
Hence, $G'\in {\cal G}_3$. As $|X|=ess'(L(\G_3))$, by (\ref{eq}), it follows that
\[
0 \ge |X|-|X'| = |X| - (|X| - |Y^3|) = |Y^3| > 0,
\]
a contradiction. This proves (ii).

Since $X$ is a minimal essential edge cut of $L(G)$ and by (i), all edges in $E(G[V_{12}(X)])$ 
receive the same color, they
are vertices in a component of $L(G)-X$, and so $G[V_{12}(G)]$ is connected subgraph of $G$.
Without loss of generality, we assume that $v_m$ is a vertex of maximum degree in $G[V_{12}(X)]$. We will show
$G[V_{12}(X)]\cong K_{1,m-1}$ to prove (iii).

We first show that $v_m$ is adjacent to each vertex $v_i$, for all $1 \le i \le m-1$. If this is false,
then there exist vertices $v_i,v_j\in V_{12}(X)$ with $i, j < m$ such that  $v_{i}v_j, v_{m}v_i\in E(G[V_{12}(X)])$
but $v_{m}v_j\notin E(G[V_{12}(X)])$.
Introduce a vertex subset $W^2$ disjoint from $V(G)$ with
$|W^2| = \max\{0, d+1-|V_1(X)|\}+1$, pick an edge $f_1 \in E_{1,X}(v_m)$ and a distinguished vertex $v'  \in W^2$. Define a new graph $G'$
(see Figure 6) as following
\begin{equation} \label{alg-9}
\begin{array}{ll}
V(G') = & V(G) \cup W^2,
\\
E(G') = & (E(G) - \{v_{i}v_{j}, f_1 \}) \cup \{v_{m}v_{j}, v_iv'\} \cup \{vw: v \in V_1(X) \mbox{ and } w \in W^2\}.
\end{array}
\end{equation}


Let  $f_2=v'v_i$, $f_3=v_iv_j$, and $f_3'=v_mv_j$. Note that $V_{12}(X)$ is also a vertex set of $G'$ and $ E_{1,X}(v_j)$,
$ E_{2,X}(v_i) - \{f_3\}$ and $E_{2,X}(v_m) - \{f_1\}$ are edge sets of $G'$. Define
\begin{equation} \label{X'-alg-9}
\begin{array}{ll}
Y_4 =  \{e_1f_3 \in X: e_1\in E_{1,X}(v_j)\}  & \mbox{ and }  Y^4 = \{e_1f_3' \in E(L(G')): e_1\in E_{1,X}(v_j)\},
\\
Y_5 = \{e_2f_3 \in X: e_2 \in E_{1,X}(v_i)\}  & \mbox{ and }  Y^5 = \{e_2'f_2 \in E(L(G'): e_2' \in E_{2,X}(v_i) - \{f_3\}\},
\\
Y_6 = \{e_3f_1 \in X: e_3 \in E_{2,X}(v_m)\}  & \mbox{ and }  Y^6 = \{e_3'f_3' \in E(L(G'): e_3' \in E_{2,X}(v_m) - \{f_1\}\}.
\end{array}
\end{equation}
Define, $X' =( X - (Y_4 \cup Y_5 \cup Y_6)) \cup Y^4 \cup Y^5 \cup Y^6$. Thus $X' = \partial_{L(G')}((E(G[V_{12}(X)]) - \{f_3\})\cup \{f_3'\})$, and so
$X'$ is a minimal essential edge cut of $L(G')$ with $V_{12}(X)=V_{12}(X')$ and $V_2(X') = \emptyset$.
By (\ref{alg-9}),  $\kappa'(G')\geq k$ and $d_{G'}(v)=d$ for each vertex $v$ in $V_{12}(X')$. Hence, $G'\in \G_3$.

\begin{center}
\begin{tikzpicture} [scale=0.5]

 \path    (0,0)  coordinate (v1);  \draw [fill=black] (v1) circle (0.07cm);    \node [right] at (v1) {\tiny $v_m$};
  \path    (0,2.5)  coordinate (vi);  \draw [fill=black] (vi) circle (0.07cm);    \node [right] at (vi) {\tiny $v_i$};
   \path    (0,5)  coordinate (vj);  \draw [fill=black] (vj) circle (0.07cm);    \node [right] at (vj) {\tiny $v_j$};

\draw (v1)--(vi)--(vj);
  \node [right] at (vj) {\tiny $v_j$};    \node [left] at (0, 3.75) {\tiny $f_3$};

\draw [dotted]  (v1)..controls (1.5,1) and  (1.5,4)..(vj);
  \begin{scope}[shift={(vi)}]
\draw (vi)--(165:1.5cm);          \draw (vi)--(195:1.5cm);
\end{scope}

  \begin{scope}[shift={(v1)}]
\draw (v1)--(165:1.5cm);          \draw (v1)--(195:1.5cm);   \draw (v1)--(180:1.5cm);
\end{scope}

  \begin{scope}[shift={(vj)}]
\draw (vj)--(165:1.5cm);          \draw (vj)--(195:1.5cm);   \draw (vj)--(180:1.5cm);
\end{scope}

 \node [above] at (4,2.5) {$\Longrightarrow$};
     \end{tikzpicture} \hskip1cm
        \begin{tikzpicture} [scale=0.5]

 \path    (0,0)  coordinate (v1);  \draw [fill=black] (v1) circle (0.07cm);    \node [right] at (v1) {\tiny $v_m$};
  \path    (0,2.5)  coordinate (vi);  \draw [fill=black] (vi) circle (0.07cm);    \node [right] at (vi) {\tiny $v_i$};
   \path    (0,5)  coordinate (vj);  \draw [fill=black] (vj) circle (0.07cm);    \node [right] at (vj) {\tiny $v_j$};

\draw (v1)--(vi);
  \node [right] at (vj) {\tiny $v_j$};    \node [left] at (0, 3.75) {\tiny $f_3$};  \node [left] at (2, 3) {\tiny $f_3'$};

\draw (v1)..controls (1.5,1) and  (1.5,4)..(vj); \draw [dotted] (vi)--(vj);

  \begin{scope}[shift={(vi)}]
\draw (vi)--(155:1.5cm);          \draw (vi)--(205:1.5cm);      \draw [blue] (vi)--(180:2cm);
 \draw [fill=black] (180:2cm) circle (0.07cm);  \node [left] at (180:2cm) {\tiny $v'$};  \node [above] at (-1.3,-0.3) {\tiny $f_2$};
\end{scope}

  \begin{scope}[shift={(v1)}]
\draw [dotted]  (v1)--(160:1.5cm);          \draw (v1)--(200:1.5cm);   \draw (v1)--(180:1.5cm);  \node [above] at (-1.3,0.3) {\tiny $f_1$};
\end{scope}

  \begin{scope}[shift={(vj)}]
\draw (vj)--(165:1.5cm);          \draw (vj)--(195:1.5cm);   \draw (vj)--(180:1.5cm);
\end{scope}

     \end{tikzpicture}

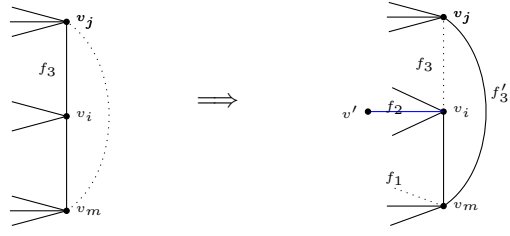
          Figure 6: Partial illustration  of forming $G'$ in (\ref{alg-9})
\end{center}

By (\ref{alg-9}), (\ref{X'-alg-9}) and the definition of $X'$, for each $v\in V_{12}(X)-\{v_m,v_i\}$,
we have $E_{1,X}(v)=E_{1,X'}(v)$ and $E_{2,X}(v)=E_{2,X'}(v)$. Moreover,  $E_{1,X'}(v_m)=E_{1,X}(v_m)-\{f_1\}$, $E_{2,X'}(v_m)=E_{2,X}(v_m) \cup \{f_3'\}$,
and $E_{1,X'}(v_i)=E_{1,X}(v_i) \cup \{f_2\}$, $E_{2,X'}(v_i)=E_{2,X}(v_i)-\{f_3\}$. By Lemma \ref{P3}(ii) and as $V_2(X') = \emptyset$, for each $v\in V_{12}(X)=V_{12}(X')$,
$|E_{1,X}(v)|+|E_{2,X}(v)|=d$ and $|E_{1,X'}(v)|+|E_{2,X'}(v)|=d$.
Since $v_m$ is a vertex of maximum degree in $G[V_{12}(X)]$,
we have $|E_{2,X}(v_i)|\le |E_{2,X}(v_m)|$ and $|E_{1,X}(v_i)| \ge |E_{1,X}(v_m))|$. Therefore, as $G' \in \G_3$ and
\begin{eqnarray*}
0 \ge |X|-|X'|
& = & \sum\limits_{v\in\{v_m,v_i\}} |E_{1,X}(v)|\cdot |E_{2,X}(v)|-\sum\limits_{v\in\{v_m,v_i\}} |E_{1,X'}(v)|\cdot |E_{2,X'}(v)| \\
&=& \big[|E_{1,X}(v_i)|\cdot |E_{2,X}(v_i)|+|E_{1,X}(v_m)|\cdot |E_{2,X}(v_m)|\big] \\
&&-\big[(|E_{1,X}(v_i)|+1)(|E_{2,X}(v_i)|-1)+(|E_{1,X}(v_m)|-1)(|E_{2,X}(v_m)|+1)\big]\\
&=& (|E_{1,X}(v_i)|-|E_{1,X}(v_m))|+(|E_{2,X}(v_m)|-|E_{2,X}(v_i)|)+2
> 0,
\end{eqnarray*}
a contradiction with $|X|=ess'(L(\G_3))$. Hence, $v_m$ must be adjacent to all other vertices in $V_{12}(X)$.

Therefore, $G[V_{12}(X)]$ contains a spanning subgraph $J \cong K_{1, m-1}$ with $V(J) = \{v_1, v_2, ..., v_m\}$ and $d_J(v_m) = m-1$.
Let $E'=E(G[V_{12}(X)])-E(J)$. If $E' = \emptyset$, then $G[V_{12}(X)]\cong K_{1,m-1}$, which proves (iii).
We are to show that if $E' \neq \emptyset$, then there will be a graph $G' \in \G_3$ that yields a contradiction to the
assumptions of the lemma.

Assume that $E' \neq \emptyset$.
Choose a set $W^3$ of vertices disjoint from $V(G)$ with $|W_3| = \max\{0, d+1-|V_1(X)|\} + 1$, a distinguished vertex $w_0 \in W^3$.
Define $G'$ to be the graph with
\begin{equation} \label{alg-10}
\begin{array}{ll}
V(G') = & V(G) \cup W^3,
\\
E(G') = & (E(G_1) - E') \cup \left(\bigcup_{v_iv_j \in E'} \{v_iw_0, v_jw_0\} \right) \cup \{vw: v \in V_1(X) \mbox{ and } w \in W^3\}.
\end{array}
\end{equation}
Hence $V_{12}(X)$ is a vertex subset of $G'$. As $J \cong K_{1,m-1}$ spans $G[V_{12}(X)]$ and by (\ref{alg-10}),
$J = G'[V_{12}(X)]$ is connected. Define $X' = \partial_{L(G')}(V(J))$. By (\ref{alg-10}),
$X'$ is a minimal essential edge cut of $L(G')$ with $V_{12}(X') = V_{12}(X)$ and $V_2(X') = \emptyset$.
Again by (\ref{alg-10}), $\kappa'(G')\geq k$ and for each vertex $v$ in $V_{12}(X')$, we have $d_{G'}(v)=d$.
Thus  $G'\in \G_3$.

\begin{center}
\begin{tikzpicture} [scale=0.5]

\draw (0,2.3) ellipse (0.6cm and 3cm);
\node [right] at (1.8,2) {\tiny $v_j$};   \node [above] at (2,3.8) {\tiny $v_i$};
\node [below] at (2,0) {\tiny $v_m$};

 \path    (2,0)  coordinate (v1);  \draw [fill=black] (v1) circle (0.07cm);  \path    (0,0)  coordinate (a1);  \draw [fill=black] (a1) circle (0.07cm);
  \path    (2,2)  coordinate (vi);  \draw [fill=black] (vi) circle (0.05cm);
    \path    (0,2)  coordinate (a21);  \draw [fill=black] (a21) circle (0.05cm); \draw (vi)--(a21);
   \path    (2,4)  coordinate (vj);  \draw [fill=black] (vj) circle (0.05cm);\draw (vi)--(vj);
         \path    (0,4)  coordinate (a3);  \draw [fill=black] (a3) circle (0.05cm);\draw (vj)--(a3);

\draw (a1)--(v1); \draw (a3)--(vj); \draw (v1)--(vj);

    \node [above] at (1,-0.2) {\tiny $e_1$};

                 \node [above] at (1,1.9) {\tiny $e_2$}; \node [above] at (1,3.9) {\tiny $e_3$};
                       \node [left] at (2.1,3) {\tiny $e_{6}$};             \node [left] at (2.1,1.1) {\tiny $e_{5}$};

\draw (vj) .. controls (3,2).. (v1);  \node [right] at (2.5,2.3) {\tiny $e_4$};

  \begin{scope}[shift={(a1)}]
\draw (a1)--(170:0.5cm);          \draw (a1)--(190:0.5cm);
\end{scope}

  \begin{scope}[shift={(a21)}]
\draw (a21)--(170:0.5cm);          \draw (a21)--(190:0.5cm);
 \end{scope}

 \begin{scope}[shift={(a3)}]
\draw (a3)--(170:0.5cm);          \draw (a3)--(190:0.5cm);
  \end{scope}

    \node [below] at (1,-1) {\tiny $V_{12}(X)=\{v_i,v_j,v_m\}$};

     \end{tikzpicture} \hskip1cm
     \begin{tikzpicture} [scale=0.5]

\draw (0,2.3) ellipse (0.6cm and 3cm);

    \path    (0,0)  coordinate (e1);  \draw [fill=black] (e1) circle (0.05cm);  \node [above] at (e1) {\tiny $e_1$};
         \path    (0,2)  coordinate (e2);  \draw [fill=black] (e2) circle (0.05cm);\node [above] at (e2) {\tiny $e_2$};
     \path    (0,3.9)  coordinate (e3);  \draw [fill=black] (e3) circle (0.05cm);\node [above] at (e3) {\tiny $e_3$};

        \path    (3,0)  coordinate (e4);  \draw [fill=black] (e4) circle (0.05cm);  \node [right] at (e4) {\tiny $e_4$};
         \path    (3,2)  coordinate (e5);  \draw [fill=black] (e5) circle (0.05cm); \node [right] at (e5) {\tiny $e_5$};
     \path    (3,3.9)  coordinate (e6);  \draw [fill=black] (e6) circle (0.05cm); \node [right] at (e6) {\tiny $e_6$};


\draw (e1)--(e4) (e1)--(e5) (e2)--(e5) (e2)--(e6) (e3)--(e4) (e3)--(e6);
\draw (e4)--(e6) (e4) .. controls (4,2)..(e6);

                                                         \begin{scope}[shift={(e1)}]
                    \draw (e1)--(170:0.5cm);          \draw (e1)--(190:0.5cm);
                                          \end{scope}

                                                                                                   \begin{scope}[shift={(e2)}]
                    \draw (e2)--(170:0.5cm);          \draw (e2)--(190:0.5cm);
                                          \end{scope}

                                                                                                   \begin{scope}[shift={(e3)}]
                    \draw (e3)--(170:0.5cm);          \draw (e3)--(190:0.5cm);
                                          \end{scope}

\draw [blue, line width=1.2] (e1)--(e4) (e1)--(e5) (e2)--(e5) (e2)--(e6) (e3)--(e6) (e3)--(e4);

    \node [below] at (1,-1) {\tiny $L(G)$};

     \end{tikzpicture} \hskip2cm
     \begin{tikzpicture} [scale=0.5]

\draw (0,2.3) ellipse (0.6cm and 3cm);
\node [right] at (1.8,2) {\tiny $v_j$};   \node [above] at (2,3.8) {\tiny $v_i$};
\node [below] at (2,0) {\tiny $v_m$};

 \path    (0,0)  coordinate (a1);  \draw [fill=black] (a1) circle (0.07cm);
     \path    (0,1.5)  coordinate (a21);  \draw [fill=black] (a21) circle (0.05cm);
     \path    (0,2.5)  coordinate (a22);  \draw [fill=black] (a22) circle (0.05cm);  \node [left] at (0.2,2.5) {\tiny $w_0$};
                   \path    (0,4)  coordinate (a3);  \draw [fill=black] (a3) circle (0.05cm);

 \path    (2,0)  coordinate (v1);  \draw [fill=black] (v1) circle (0.07cm);
  \path    (2,2)  coordinate (vi);  \draw [fill=black] (vi) circle (0.05cm);
   \path    (2,4)  coordinate (vj);  \draw [fill=black] (vj) circle (0.05cm);
\draw (a22)--(vj) (a22)--(vi); \draw (a21)--(vi);
\draw (vj)--(a3);
\draw (a1)--(v1); \draw (a3)--(vj); \draw (v1)--(vi);
\draw [dotted] (vi)--(vj);

    \node [above] at (1,-0.2) {\tiny $e_1$};

                 \node [above] at (1,1.2) {\tiny $e_2$}; \node [above] at (1,3.9) {\tiny $e_3$};
                           \node [left] at (2.1,1.1) {\tiny $e_{5}$};

\draw (vj) .. controls (3,2).. (v1);  \node [right] at (2.5,2.3) {\tiny $e_4$};

  \begin{scope}[shift={(a1)}]
\draw (a1)--(170:0.5cm);          \draw (a1)--(190:0.5cm);
\end{scope}

  \begin{scope}[shift={(a21)}]
\draw (a21)--(170:0.5cm);          \draw (a21)--(190:0.5cm);
 \end{scope}

 \begin{scope}[shift={(a3)}]
\draw (a3)--(170:0.5cm);          \draw (a3)--(190:0.5cm);
  \end{scope}

    \node [below] at (1,-1) {\tiny $G'$};

     \end{tikzpicture} \hskip1cm
     \begin{tikzpicture} [scale=0.5]

\draw (0,2.3) ellipse (0.6cm and 3cm);

    \path    (0,0)  coordinate (e1);  \draw [fill=black] (e1) circle (0.05cm);  \node [above] at (e1) {\tiny $e_1$};
         \path    (0,1)  coordinate (e2);  \draw [fill=black] (e2) circle (0.05cm);\node [above] at (e2) {\tiny $e_2$};
     \path    (0,2)  coordinate (e-1);  \draw [fill=black] (e-1) circle (0.05cm);\node [left] at (e-1) {\tiny $w_0v_j$};
\path    (0, 3)  coordinate (e-2);  \draw [fill=black] (e-2) circle (0.05cm);  \node [left] at (e-2) {\tiny $w_0v_i$};
\path    (0,4)  coordinate (e3);  \draw [fill=black] (e3) circle (0.05cm);  \node [right] at (e3) {\tiny $e_3$};

         \path    (3,1)  coordinate (e5);  \draw [fill=black] (e5) circle (0.05cm); \node [right] at (e5) {\tiny $e_5$};
     \path    (3,3)  coordinate (e4);  \draw [fill=black] (e4) circle (0.05cm); \node [right] at (e4) {\tiny $e_4$};

\draw (e1)--(e4) (e1)--(e5) (e2)--(e5) (e-1)--(e5) (e-2)--(e3) (e-2)--(e4) (e3)--(e4) (e4)--(e5);

                                                         \begin{scope}[shift={(e1)}]
                    \draw (e1)--(170:0.5cm);          \draw (e1)--(190:0.5cm);
                                          \end{scope}

                                                                                                   \begin{scope}[shift={(e2)}]
                    \draw (e2)--(170:0.5cm);          \draw (e2)--(190:0.5cm);
                                          \end{scope}

                                                                                                   \begin{scope}[shift={(e3)}]
                    \draw (e3)--(170:0.5cm);          \draw (e3)--(190:0.5cm);
                                          \end{scope}


    \node [below] at (1,-1) {\tiny $L(G')$};

     \end{tikzpicture}

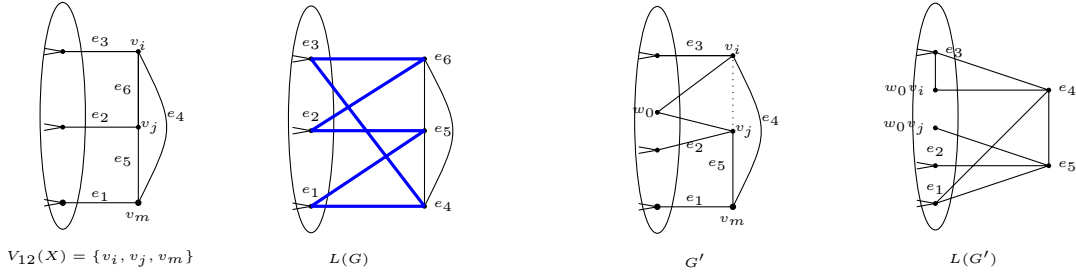
     Figure 7: Illustration of forming $G'$ in (\ref{alg-10}) with $X=\{e_1e_4,e_1e_5,e_2e_5,e_2e_6,e_3e_4,e_3e_6\}$.
\end{center}

(See Figure 7). By (\ref{alg-10}), for $j \in \{1,2\}$,  $E_{j,X'}(v_m)=E_{j,X}(v_m)$ in $G'$.
As $G'[V_{12}(X')] \cong K_{1, m-1}$,  for each $v\in V_{12}(X')-\{v_m\}$, we have $|E_{1,X'}(v)|=d-1$ and $|E_{2,X'}(v)|=1$.
By Observation \ref{Obs(G)}(iii), for each $v\in \{v_1,v_2,\cdots,v_{m-1}\}$, $|E_{1,X}(v)|+|E_{2,X}(v)|=d$. Since
$ (|E_{1,X}(v)|-1) (|E_{2,X}(v)|-1)\ge 0$, it follows that  $|E_{1,X}(v)|\cdot |E_{2,X}(v)|\ge d-1$. Thus,
\begin{eqnarray*}
 |X|-|X'|
& = & \sum\limits_{v\in V_{12}(X)} |E_{1,X}(v)|\cdot |E_{2,X}(v)|-\sum\limits_{v\in V_{12}(X')} |E_{1,X'}(v)|\cdot |E_{2,X'}(v)| \\
&=& \sum\limits_{v\in \{v_1,v_2,\cdots,v_{m-1}\}} |E_{1,X}(v)|\cdot |E_{2,X}(v)|- (m-1) (d-1) \\
&\ge & (m-1)(d-1)-(m-1)(d-1)=0
\end{eqnarray*}
Since $|X|=ess'(L({\cal G}_3))$, we have $|X|=|X'|$. Therefore,  the set $X'$ forms an $ess'(G')$-edge cut of $L(G')$ with $V_{12}(X')\cong K_{1,m-1}$,
contrary to the hypothesis that
$|E(G[V_{12}(X)])|$ is minimized.   Hence, $G[V_{12}(X)]\cong K_{1,m-1}$. This completes the proof of (iii), as well as the lemma.
\q

Let $j\in\{1,2\}$. Define
$${\cal G}_j'=\left\{\begin{array}{l|l}
 G\in {\cal G}_j  & G \hbox{  has a minimum essential edge-cut of  }
 X=X(G)  \hbox{ with } |X|=ess'(L(\G_j)),  \\
  &  \hbox{ and  assume, subject to this, that } |E(G[V_{12}(X)])|
 \hbox{  is minimized. }
\end{array}\right\}, $$
and define
$${\cal G}'_3=\left\{\begin{array}{l|l}
 G\in {\cal G}_3\cap {\cal G}_0(d,k)  & G \hbox{  has a minimum essential edge-cut of  }
 X=X(G)  \hbox{ with } |X|=ess'(L(\G_3)), \\
  &  \hbox{ and  assume, subject to this, that } |E(G[V_{12}(X)])|  \hbox{  is minimized. }
\end{array}\right\}.$$

By Observation \ref{Obs(G)}(iv), $\G_1\cup\G_2\cup\G_3=\G$. As $\G_1$, $\G_2$ and $\G_3$ are
pairwise disjoint, we conclude that
\begin{equation} \label{min}
\kappa'_{L^2}(d,k) = \min\{ess'(L(\G'_i)): 1 \le i \le 3\}.
\end{equation}


\noindent
{\bf{Proof of Theorem \ref{main-1}.}}   Let $G\in \G_1'\cup\G_2'\cup\G_3'$.
Let $X$ be an $ess'(L(G))$-essential edge cut of $L(G)$ with $m = |V_{12}(X)|$ in $G$, where $1\leq m\leq k$. Let $V(G[V_{12}(X)])=\{v_1,v_2,\cdots,v_m\}$.

Suppose that $G\in\G'_3$.
Then $m\geq 3$, and otherwise $G$ has no essential edge-cut.
Without loss of generality, assume that $V_1(X)\neq\emptyset$ and $V_2(X)=\emptyset$.
By Lemma \ref{P3}, $G[V_{12}(X)]$
is isomorphic to $K_{1,m-1}$ and $d_G(v_i)=d$ for $1\leq i\leq m$, with $d_{G[V_{12}(X)]}(v_i)=1$, $1 \le i \le m-1$,  and $d_{G[V_{12}(X)]}(v_{m})=m-1$.
Thus $|E_{2,X}(v_i)|=1$ for $i=1,2,\cdots,m-1$ and $|E_{2,X}(v_m)|=m-1$.
By \eqref{eq},
\begin{equation}
\begin{split}
|X|&=\sum \limits_{v\in V_{12}(X)}|E_{1,X}(v)|\cdot |E_{2,X}(v)|=|E_{1,X}(v_{m})| \cdot|E_{2,X}(v_{m})|+\sum \limits_{i=1}^{m-1} |E_{1,X}(v)|\cdot |E_{2,X}(v)| \\
&=(d-(m-1))(m-1)+\sum \limits_{i=1}^{m-1} |E_1(v)|=(d-m+1)(m-1)+ (m-1)(d-1)\\
&=-m^2+(2d+1)m-2d.\nonumber
\end{split}
\end{equation}
Let $f(m)=-m^2+(2d+1)m-2d$.
As $d\ge k$ and $d\ge 3$, we have  $2d-2\geq k$, and so $3\le m\le k\le 2d-2$.
Hence $f(m)$ on $[3,2d-2]$ has a minimum value $f(3)=f(2d-2)=4d-6$, and so $ess'(L(\G_3)) = |X|=f(m)\ge 4d-6$.

Suppose that $V_1(X)\cup V_2(X)\neq \emptyset$.
If
$G\in {\cal G}_1'$, then  there exists $i_0\in\{1,2,\cdots,m\}$ such that $a_{i_0,1,X}=a_{i_0,2,X}=\lceil\frac{d}{2}\rceil$.
By Proposition \ref{edge cut} and Lemma \ref{P1}(i) and (ii),
$\sum\limits_{i=1}^m a_{i,1,X}=k$.
By \eqref{eq}, we have
\begin{equation}\label{G1=}
\begin{split}
|X|
=&\sum \limits_{v\in V_{12}(X)}|E_{1,X}(v)|\cdot |E_{2,X}(v)|=\Big(\left\lceil \frac{d}{2} \right\rceil\Big)^2 +
\sum \limits_{1\leq i\neq i_0\leq m}a_{i,1,X}(d-a_{i,1,X})
\\
=&\Big(\left\lceil \frac{d}{2} \right\rceil\Big)^2 +
d\Big(\sum \limits_{i=1}^m   a_{i,1,X}-a_{i_0,1,X}\Big)-\sum \limits_{1\leq i\neq i_0\leq m}a_{i,1,X}^2
\\
=&\Big(\left\lceil \frac{d}{2} \right\rceil\Big)^2
+\Big(k-\left\lceil \frac{d}{2} \right\rceil\Big)d-
\sum \limits_{1\leq i\neq i_0\leq m}a^{2}_{i,1,X}.
\end{split}
\end{equation}

If  $G\in \G'_2$, by Lemma \ref{P1}(v), $d_G(v_i)=a_{i,1,X}+a_{i,2,X}=d$.  By \eqref{eq},  we have
\begin{equation}\label{G2=}
\begin{split}
|X|&=\sum \limits_{v\in V_{12}(X)} |E_{1,X}(v)|\cdot |E_{2,X}(v)|
=\sum \limits_{i=1}\limits^{m}a_{i,1,X}(d-a_{i,1,X})
=kd-\sum \limits_{i=1}\limits^{m}a^{2}_{i,1,X}.
\end{split}
\end{equation}

\noindent
{\bf Case 1. }  $1\leq k\leq \lfloor\frac{d}{2}\rfloor$.

Suppose $G\in\G'_1$.  As $a_{i_0,1,X}=a_{i_0,2,X}=\lceil\frac{d}{2}\rceil$ and $\sum\limits_{i=1}^m a_{i,1,X}=k$,  we have $k\ge a_{i_0,1,X}=\lceil\frac{d}{2}\rceil$. Thus
$k=\frac{d}{2}$ and so $d$ is even. As $\frac{d}{2}=a_{x_0,1,X}\le \sum\limits_{i=1}^m a_{i,1,X}=k=\frac{d}{2}$, we have $a_{i,1,X}=0$ for $1\leq i\neq i_0\leq m$. By \eqref{G1=},  we have
\begin{equation}
|X|=\Big(\left\lceil \frac{d}{2} \right\rceil\Big)^2
+\Big(k-\left\lceil \frac{d}{2} \right\rceil\Big)d-
\sum \limits_{1\leq i\neq i_0\leq m}a^{2}_{i,1,X}
=\frac{d^2}{4}=f(d,k). \nonumber
\end{equation}

Now we suppose $G\in \G'_2$.  Then $\sum\limits^{m}\limits_{i=1}a_{i,1,X}=k$.
 By \eqref{G2=}, we have
\begin{equation}
\begin{split}
|X|=kd-\sum \limits_{i=1}\limits^{m}a^{2}_{i,1,X} \ge kd- (\sum \limits_{i=1}\limits^{m}a_{i,1,X})^2=kd-k^2=f(d,k).\nonumber
\end{split}
\end{equation}
So,  if $V_1(G),V_2(G)\neq \emptyset$
and $1\leq k\leq \lfloor\frac{d}{2}\rfloor$, then $ess'({L(\G_1)})=f(d,k)$
and $ess'({L(\G_2)})\geq f(d,k)$.

\noindent
{\bf Case 2.  } $\lfloor\frac{d}{2}\rfloor< k < d$.

Then $k\leq 2\lfloor\frac{d}{2}\rfloor$.
Suppose $G\in \G_1'$, by \eqref{G1=} and
$\sum\limits_{1\leq i\neq i_0\leq m}a_{i,1,X}=k-\lceil \frac{d}{2} \rceil$, we have
\begin{equation}
\begin{split}
|X|
&=\Big(\left\lceil \frac{d}{2} \right\rceil\Big)^2
+\Big(k-\left\lceil \frac{d}{2} \right\rceil\Big)d-
\sum \limits_{1\leq i\neq x_0\leq m} a^{2}_{i,1,X}\\
&\geq \Big(\left\lceil \frac{d}{2} \right\rceil\Big)
^2+\Big(k-\left\lceil \frac{d}{2}\right\rceil\Big)d-
\Big(k-\left\lceil \frac{d}{2} \right\rceil\Big)^2
=kd-k^2+2k\Big(\left\lceil \frac{d}{2}\right\rceil\Big)-
\Big(\left\lceil\frac{d}{2}\right\rceil\Big)d. \nonumber
\end{split}
\end{equation}
Now, we suppose $G\in \G_2'$.  As $d_G(v_i)=a_{i,1,X}+a_{i,2,X}=d$ and $a_{i,1,X}\le a_{i,2,X}$,
we have $a_{i,1,X}\le\lfloor \frac{d}{2}\rfloor$. By Lemma \ref{P1}(ii), $\sum\limits_{i=1}^m a_{i,1,X}=k$.
Since the function $x^2$ is convex, by the majorization principle (see \cite{HLP}),  the maximum of $\sum\limits_{i=1}^m a_{i,1,X}^2$
under the given constraints $\sum\limits_{i=1}^m a_{i,1,X}=k$ occurs when as many variables as possible take the maximal value
$\lfloor\frac{d}{2}\rfloor$ and at most one variable takes the remaining value.

Let $t=\lfloor \frac{d}{2}\rfloor$. Write $k=qt+r$, where $q=\lfloor \frac{k}{t}\rfloor$ and $0\le r< t$. The maximal possible $\sum\limits_{i=1}^m a_{i,1,X}^2$ is
$$S_{max}=qt^2+r^2.$$
Notice that $k^2-2kt+2t^2=(k-t)^2+t^2$. Substitute $k=qt+r$, then
$$(k-t)^2+t^2=(q-1)^2 t^2+2(q-1) tr+r^2+t^2.$$
Therefore,
$$(k-t)^2+t^2-S_{max}=(q-1) t[(q-2)t+2r].$$
Since $q\ge 2$ and $r\ge 0$, $(q-1) t[(q-2)t+2r]\ge 0$. Thus, $S_{max}\le (k-t)^2+t^2$.  By \eqref{G2=}, we have
$$
|X|=kd-\sum\limits_{i=1}^m a_{i,1,X}^2 \ge kd-\Big[\Big(k-\Big\lfloor \frac{d}{2}\Big\rfloor\Big)^2+\Big\lfloor \frac{d}{2}\Big\rfloor^2\Big]=
 kd-k^2+2k\left\lfloor \frac{d}{2}\right\rfloor-2\Big(\left\lfloor \frac{d}{2}\right\rfloor\Big)^2.
$$

That is $ess'({L(\G_1)})\geq f(d,k)$, and $ess'({L(\G_2)})\geq f(d,k)$.

Moreover, if $\lfloor\frac{d}{2}\rfloor< k <\frac{3d-1}{4}$, 
$kd-k^2+2k\left\lceil \frac{d}{2}\right\rceil-
\left\lceil\frac{d}{2}\right\rceil d\leq kd-k^2+2k\left\lfloor \frac{d}{2}\right\rfloor-2\left\lfloor \frac{d}{2}\right\rfloor^2 $, then
$ess'({L(\G_1)})\leq ess'({L(\G_2)})$;
if $\frac{3d-1}{4}\leq k<d$, $ess'({L(\G_1)}) \geq ess'({L(\G_2)})$.

\noindent
{\bf Case 3. } $d=k$.

Suppose $G\in \G_1'$. By \eqref{G1=} and
$\sum\limits_{1\leq i\neq x_0\leq m}a_{i,1,X}=k-\lceil \frac{d}{2} \rceil=\lfloor \frac{d}{2}\rfloor$, we have
\begin{equation}
\begin{split}
|X|
&=\Big(\left\lceil \frac{d}{2} \right\rceil\Big)^2
+\Big(k-\left\lceil \frac{d}{2} \right\rceil\Big)d-
\sum \limits_{1\leq i\neq x_0\leq m}a^{2}_{i,1,X}\\
&=\Big(\left\lceil \frac{d}{2} \right\rceil\Big)^2
+\Big(\left\lfloor \frac{d}{2} \right\rfloor\Big)d-
\sum \limits_{1\leq i\neq x_0\leq m}a^{2}_{i,1,X}\\
&\geq \Big(\left\lceil \frac{d}{2} \right\rceil\Big)
^2+\Big(\left\lfloor \frac{d}{2} \right\rfloor\Big)d-
\Big(\left\lfloor \frac{d}{2} \right\rfloor\Big)^2
=d\Big(\left\lceil \frac{d}{2} \right\rceil\Big)=f(d,k). \nonumber
\end{split}
\end{equation}

Now, suppose that $G\in \G_2'$. By \eqref{G2=} and
$\sum\limits_{i=1}\limits^{m}a_{i,1,X}=k=d$ with $a_{i,1.X}\leq \lfloor \frac{d}{2} \rfloor$. Let $b= \lfloor \frac{d}{2} \rfloor$.
Since the function $x^2$ is convex, by the majorization principle (see \cite{HLP}),  the sum $\sum\limits_{i=1}^m a_{i,1,X}^2$
with  $\sum\limits_{i=1}^m a_{i,1,X}=d$ and $0<a_{i,1,X}\le b$ is maximized by taking the the largest possible number of parts equal to
$b$  and one remaining part equal to the remainder. As indicated in \cite{HLP}, expressing $d=2b+r$ with $r\in\{0,1\}$, the maximum is
$$\max \left\{\sum\limits_{i=1}^m a_{i,1,X}^2 \right\}=2b^2+r^2.$$
In particular this yields $2b^2$ when $d$ is even and $2b^2+1$ when $d$ is odd.

By \eqref{G2=}, we have
$$
|X|=d^2-\sum\limits_{i=1}^m a_{i,1,X}^2\ge
\left\{\begin{array}{ll}
d^2-2\lfloor\frac{d}{2}\rfloor^2=\frac{1}{2}d^2=d\cdot\lfloor \frac{d}{2}\rfloor, & \hbox{ if } d \hbox{ is even}, \\
d^2-(2\lfloor\frac{d}{2}\rfloor^2+1)=\frac{d^2+2d-3}{2}=\frac{d(d+1)+(d-3)}{2} \ge d\cdot \lceil \frac{d}{2} \rceil, & \hbox{ if } d \hbox{ is odd}.
\end{array}\right.
$$
Therefore,  $ess'(L(\G_1))\le ess'(L(\G_2))$.

As shown in the proofs above,  $\kappa_{L^2}(d,k)\geq min\{f(d,k),4d-6\}$.
We now construct examples to demonstrate
that $\kappa_{L^2}(d,k)\leq \min\{f(d,k),4d-6\}$.

Consider the graph  $J(d,k,n_1,n_2)$ in Example \ref{ex-1}. In this case,  both $V_1(G)$ and $V_2(G)$ are non-empty.
If $1\leq k\leq \lfloor\frac{d}{2}\rfloor$,  then
by (2),  $\kappa_{L^2}(d,k)\leq |X|=k(d-k)$.

When $\lfloor\frac{d}{2}\rfloor< k\le d$, (see Figure 8),  let $v_1, v_2\in V_{12}(X)$, where $G[V_1(X)]$, $G[V_2(X)]$ are two
complete subgraphs of order $d+1$ induced by $G$. If $\lfloor\frac{d}{2}\rfloor< k <\frac{3d-1}{4}$, then $v_1$ is adjacent to $\lceil\frac{d}{2}\rceil$ vertices in $V_1(X)$ and to  $\lceil\frac{d}{2}\rceil$ vertices in $V_2(X)$;
while  $v_2$ is adjacent to $k-\lceil\frac{d}{2}\rceil$
vertices in $V_1(X)$ and $d-k+\lceil\frac{d}{2}\rceil$ vertices in $V_2(X)$. If $\frac{3d-1}{4}\leq k<d$, then  $v_1$ is adjacent to $\lfloor\frac{d}{2}\rfloor$ vertices in $V_1(X)$
and to $\lfloor\frac{d}{2}\rfloor$ vertices in $V_2(X)$, while  $v_2$ is adjacent to $k-\lfloor\frac{d}{2}\rfloor$
vertices in $V_1(X)$ and $d-k+\lfloor\frac{d}{2}\rfloor$ vertices in $V_2(X)$. If $k=d$, then $v_1$ is adjacent to $\lceil\frac{d}{2}\rceil$ vertices in $V_1(X)$
and $\lceil\frac{d}{2}\rceil$ vertices in $V_2(X)$, while
$v_2$ is adjacent to $d-\lceil\frac{d}{2}\rceil$
vertices in $V_1(G)$ and to $\lceil\frac{d}{2}\rceil$ vertices in $V_2(X)$.

\begin{center}
       \begin{tikzpicture} [scale=0.5]

   \path    (3,1.1)  coordinate (a);  \draw [fill=black] (a) circle (0.08cm); \node [below] at (a) { $v_1$};
      \path    (3, -1.1)  coordinate (b);  \draw [fill=black] (b) circle (0.08cm); \node [below] at (b) { $v_2$};
 \draw (0,0) ellipse (1.5cm and 3cm);       \draw (6,0) ellipse (1.5cm and 3cm);

  \path    (0,0.5)  coordinate (v1);  \draw [fill=black] (v1) circle (0.08cm); \node at (0,1.5) {$\vdots$};
    \path    (0,2)  coordinate (v2);  \draw [fill=black] (v2) circle (0.08cm);
    \draw (a)--(v1) (a)--(v2);

    \path    (0,-0.5)  coordinate (v3);  \draw [fill=black] (v3) circle (0.08cm); \node at (0,-1) {$\vdots$};
    \path    (0,-2)  coordinate (v4);  \draw [fill=black] (v4) circle (0.08cm);
     \draw (b)--(v3) (b)--(v4);

  \path    (6,0.5)  coordinate (u1);  \draw [fill=black] (u1) circle (0.08cm); \node at (6,1.5) {$\vdots$};
    \path    (6,2)  coordinate (u2);  \draw [fill=black] (u2) circle (0.08cm);
    \draw (a)--(u1) (a)--(u2);

    \path    (6,-0.5)  coordinate (u3);  \draw [fill=black] (u3) circle (0.08cm); \node at (6,-1) {$\vdots$};
    \path    (6,-2)  coordinate (u4);  \draw [fill=black] (u4) circle (0.08cm);
     \draw (b)--(u3) (b)--(u4);

\node [below] at (0,-3.5) {$V_1(G)$}; \node [below] at (6,-3.5) {$V_2(G)$};

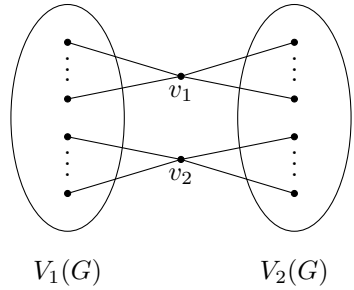
 \node [below] at (3,-5) {Figure 8.  Illustration of $\lfloor \frac{d}{2}\rfloor <k\le d$};
     \end{tikzpicture}

\end{center}

Therefore, by (\ref{eq}),  $$
\kappa_{L^2}(d,k)\leq \left\{\begin{array}{ll}
kd-k^2+2k(\lceil \frac{d}{2}\rceil)-(\lceil\frac{d}{2}\rceil)d,  & \hbox{ if } \lfloor\frac{d}{2}\rfloor< k <\frac{3d-1}{4}, \\
kd-k^2+2k\lfloor \frac{d}{2}\rfloor-(\lfloor \frac{d}{2}\rfloor) d \le kd-k^2+2k\lfloor \frac{d}{2}\rfloor-2(\lfloor \frac{d}{2}\rfloor)^2, & \hbox{ if } \frac{3d-1}{4}\leq k<d, \\
d(\lceil \frac{d}{2}\rceil), & \hbox{ if } k=d.
\end{array}\right.$$

\begin{center}
       \begin{tikzpicture} [scale=0.45]

     \path    (5,3)  coordinate (c);  \draw [fill=black] (c) circle (0.1cm); \node [right] at (c) { $v_1$};
   \path    (5,0.5)  coordinate (a);  \draw [fill=black] (a) circle (0.1cm); \node [right] at (a) { $v_2$};
      \path    (5, -2)  coordinate (b);  \draw [fill=black] (b) circle (0.1cm); \node [right] at (b) { $v_3$};
 \draw (0,0.2) ellipse (1.8cm and 4cm);

   \path    (0, 3.5)  coordinate (v5);  \draw [fill=black] (v5) circle (0.1cm); \node at (0,3) {$\vdots$};
    \path    (0,2)  coordinate (v6);  \draw [fill=black] (v6) circle (0.1cm);

     \path    (0,1)  coordinate (v2);  \draw [fill=black] (v2) circle (0.1cm);
  \path    (0,-0.5)  coordinate (v1);  \draw [fill=black] (v1) circle (0.1cm); \node at (0,0.5) {$\vdots$};

    \path    (0,-1.5)  coordinate (v3);  \draw [fill=black] (v3) circle (0.1cm); \node at (0,-2) {$\vdots$};
    \path    (0,-3)  coordinate (v4);  \draw [fill=black] (v4) circle (0.1cm);

      \draw (a)--(v1) (a)--(v2);
     \draw (b)--(v3) (b)--(v4);
          \draw (c)--(v5) (c)--(v6);
     \draw (b)--(c);

\node [below] at (0,-4) {$V_1(G)$};

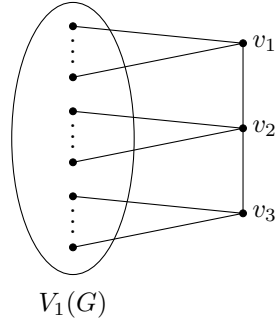
 \node [below] at (3,-6) {Figure 9:   Illustration of $V_1(X)\not=\emptyset$ and  $V_2(X)=\emptyset$};
     \end{tikzpicture}

\end{center}

Finally, when $V_1(X)\neq \emptyset$ and $V_2(X)=\emptyset$,
let $v_1, v_2, v_3\in V_{12}(X)$,  and let $G[V_1(X)]$ be a
complete graphs of order $d+1$ induced by $G$. Suppose
$G[\{v_1,v_2,v_3\}]$ is a path, and
the vertices $v_1$ and $v_3$ are adjacent to $d-1$
vertices in $V_1(X)$, and
$v_2$ is adjacent with $d-2$ vertex in $V_1(X)$, (see Figure 9).  By \eqref{eq},  $|X|\leq(d-1)+2(d-2)+(d-1)=4d-6$, and hence $\kappa_{L^2}(d,k)\leq 4d-6$.
\q

\section{Concluding Remark}
Line graphs have been considered as commonly adopted interconnection network models. As various kinds of graph connectivity reflect the reliabilities of networks,
the investigation of iterated line graph connectivity has become one of the attractive research problems. The current research extends
the former result on $\kappa_{L^2}(d,1)$ of Knor and Niepel \cite{KnNi03} to a general form of $\kappa_{L^2}(d,k)$ for all integers $d \ge 3$ and $k \ge 1$.
Motivated by Shao's recent sequence of work in \cite{Shao05, Shao10, Shao18, Shao22+}, performing further study may be needed to understand the growth trend or pattern of the connectivity of
iterated line graphs. Formally, for any integer triples $(d, k, n)$ with $d \ge 3$, $k \ge 1$ and $n \ge 1$, define $\G(d,k)$ as in (\ref{def-g-k}) and define a function
\[
\phi(d,k, n)  = \inf\{ \kappa(L^n(G)): G\in \G(d,k)\}.
\]
It is of interests to know, when $d$ and $k$ are fixed, how the function $\phi(d,k, n)$ grows as $n$ grows, and how the variables $d$ and $k$ impact the behavior of
the function $\phi(d,k, n)$. All these remain to be investigated in the future.

\noindent
\textbf{Acknowledgements}
Xiong was supported by the National Natural Science Foundation of China (No.12461067). Lai was supported by the National Natural Science Foundation of China (No.12471333).

\section*{Declarations}

\noindent
\textbf{Conflict of interest}
The author has no relevant financial or non-financial interests to disclose.

\end{document}